\documentclass[11pt, reqno]{amsart}
\usepackage{amsfonts,latexsym}
\usepackage{amsmath}
\usepackage{amscd}
\usepackage{float,amsmath,amssymb,mathrsfs,bm,multirow,graphics}
\usepackage[dvips]{graphicx}
\usepackage[percent]{overpic}

\newcommand{\nequation}{\setcounter{equation}{0}}
\renewcommand{\theequation}{\mbox{\arabic{section}.\arabic{equation}}}
\newcommand{\R}{{\Bbb R}}

\newcommand{\C}{{\Bbb C}}

\newcommand{\proofbegin}{\noindent{\it Proof.\,\,}}
\newcommand{\proofend}{\hfill$\Box$\bigskip}

\newcommand{\lot}{\text{\upshape lower order terms}}
\newcommand{\re}{\text{\upshape Re\,}}
\newcommand{\res}{\text{\upshape Res\,}}

\def\Xint#1{\mathchoice
{\XXint\displaystyle\textstyle{#1}}%
{\XXint\textstyle\scriptstyle{#1}}%
{\XXint\scriptstyle\scriptscriptstyle{#1}}%
{\XXint\scriptscriptstyle\scriptscriptstyle{#1}}%
\!\int}
\def\XXint#1#2#3{{\setbox0=\hbox{$#1{#2#3}{\int}$}
\vcenter{\hbox{$#2#3$}}\kern-.5\wd0}}

\def\dashint{\Xint-}

\allowdisplaybreaks

\newtheorem{theorem}{Theorem}[section]

\newtheorem{lemma}[theorem]{Lemma}

\newtheorem{remark}[theorem]{Remark}
\newtheorem{example}[theorem]{Example}
\newtheorem{figuretext}{Figure}

\input epsf
\title[The Unified Method I]{\sc The Unified Method: I~Non-Linearizable Problems on the Half-Line}
\author{A. S. Fokas}
\address{A.S.F.: Department of Applied Mathematics and Theoretical Physics, University of Cambridge, Cambridge CB3 0WA, United Kingdom, and Research Center of Mathematics, Academy of Athens, 11527, Greece.}
\email{T.Fokas@damtp.cam.ac.uk} 

\author{J. Lenells}
\address{J.L.: Department of Mathematics, Baylor University, One Bear Place \#97328, Waco, TX 76798, USA.}
\email{Jonatan\_Lenells@baylor.edu}

\begin{document}

\begin{abstract} 
\noindent
Boundary value problems for integrable nonlinear evolution PDEs formulated on the half-line can be analyzed by the unified method introduced by one of the authors and used extensively in the literature. The implementation of this general method to this particular class of problems yields the solution in terms of the unique solution of a matrix Riemann-Hilbert problem formulated in the complex $k$-plane (the Fourier plane), which has a jump matrix with explicit $(x,t)$-dependence involving four scalar functions of $k$, called spectral functions. Two of these functions depend on the initial data, whereas the other two depend on all boundary values. The most difficult step of the new method is the characterization of the latter two spectral functions in terms of the given initial and boundary data, i.e. the elimination of the unknown boundary values. For certain boundary conditions, called linearizable, this can be achieved simply using algebraic manipulations. Here, we present an effective characterization of the spectral functions in terms of the given initial and boundary data for the general case of non-linearizable boundary conditions. This characterization is based on the analysis of the so-called global relation, on the analysis of the equations obtained from the global relation via certain transformations leaving the dispersion relation of the associated linearized PDE invariant, and on the computation of the large $k$ asymptotics of the eigenfunctions defining the relevant spectral functions.
\end{abstract}

\maketitle




\section{Introduction}\nequation
A unified method for analyzing boundary value problems, extending ideas of the so-called inverse scattering transform method, was introduced in \cite{F1997, F2000, Fbook}. This method has been implemented to linear and integrable nonlinear evolutions PDEs on the half-line and the finite interval \cite{
BF2008, BFS2004, BFS2006, BS2003, D2007, DTVpreprint, D2009, F2002, F2002IMA, F2004, F2007, FG1994, FI2004, FIS2005, FL2010, FP2005, I2003, K2010, Ldnls, LdnlsD2N, LFgnls, P2004, P2005, P2006, P2007, S2005, Spreprint, V2007} to linear and integrable nonlinear evolution PDEs in two space variables \cite{F2002IMA, F2003, F2009, KF2010, MF2011}, to linear elliptic PDEs 
\cite{AF2005, bAF1999, bAF2001, DF2005, D2010, FK2000, FK2003, FP2006, HZ2011, SF2010b}, and to the two prototypical integrable nonlinear two-dimensional elliptic PDEs, namely the sine-Gordon and Ernst equations \cite{FPpreprint, FPLpreprint, Lholedisk, LFernst, PP2010}.

Regarding {\it linear} evolution PDEs containing spatial derivatives of {\it second} order, the unified method yields novel integral representations which have both analytic and numerical advantages in comparison with the classical integral and series representations: (i) They are uniformly convergent at the boundaries; analytically, this makes it easier to prove rigorously the validity of such representations without the a priori assumption of existence, whereas numerically, using appropriate contour deformations, it makes it possible to obtain integrands which decay exponentially as $|k|\to \infty$ and this leads to efficient numerical computations \cite{FF2008, PK2009}. (ii) These representations retain their form even for more complicated boundary conditions, whereas the classical representations involve infinite series over a spectrum determined by a transcendental equation. For example, in the case of the heat equation with Robin boundary conditions, the classical representation involves an infinite {\it series} over $\{k_n\}_0^\infty$ where $k_n$ satisfy a transcendental equation, whereas the new method yields an explicit \textit{integral} representation. For linear evolution PDEs containing $x$-derivatives of {\it arbitrary} order formulated on the half-line or the finite interval, the unified method still yields integral representations explicitly defined in terms of the given initial and boundary conditions, whereas it is shown in \cite{P2011} and \cite{S2011} that even for the linearized KdV on the finite interval with generic boundary conditions, there does \textit{not} exist an infinite series representation.

Regarding {\it integrable nonlinear} evolution PDEs in one space variable, the unified method yields novel integral representations formulated in the complex $k$-plane (the Fourier plane). These representations are similar to the integral representations for the linearized versions of these nonlinear PDEs, but also contain the entries of a certain matrix-valued function $M(x,t,k)$, which is the solution of a matrix Riemann-Hilbert (RH) problem. The main advantage of the new method is the fact that this RH problem involves a jump matrix with {\it explicit} $(x,t)$-dependence, uniquely defined in terms of four scalar functions called {\it spectral functions} and denoted by $\{a(k), b(k), A(k), B(k)\}$. The functions $a(k)$ and $b(k)$ are defined in terms of the initial data $q_0(x) = q(x,0)$ via a system of linear Volterra integral equations. The functions $A(k)$ and $B(k)$ are also defined via a system of linear Volterra integral equations, but these integral equations involve {\it all} boundary values. For example, for the nonlinear Schr\"odinger (NLS) and the modified Korteweg-de Vries (mKdV) equations formulated on the half-line, $A(k)$ and $B(k)$ are defined in terms of $\{q(0,t), q_x(0,t)\}$ and $\{q(0,t), q_x(0,t), q_{xx}(0,t)\}$ respectively. A major difficulty of initial-boundary value problems is that some of these boundary values are unknown. For example, for the Dirichlet problem of the NLS and of the mKdV, the boundary values $q_x(0,t)$ and $\{q_x(0,t), q_{xx}(0,t)\}$ respectively, are unknown. It turns out that this difficulty can be bypassed by utilising the so-called \textit{global relation}, which is a simple algebraic equation which couples the spectral functions.  

\subsection{Linearizable versus non-linearizable BVPs}
We now distinguish two cases: $(a)$ For a particular class of boundary conditions called {\it linearizable}, it is possible to express $A(k)$ and $B(k)$ in terms of $\{a(k), b(k)\}$ and the given boundary conditions. This can be achieved by analyzing the global relation and the equations obtained from the global relation under those transformations in the complex $k$-plane which leave invariant the dispersion relation of the linearized version of the given nonlinear PDE. It must be emphasized that the above analysis is carried out in the $k$-space, thus the spectral functions $A(k)$ and $B(k)$ are determined {\it directly}, without the need of determining the unknown boundary values. In summary, for this class of problems the unified method is as effective as the classical inverse scattering transform method. $(b)$ For non-linearizable boundary conditions which decay for large $t$, by utilizing the crucial feature of the new method that it yields RH problems with explicit $(x,t)$-dependence, it is possible to obtain useful asymptotic information about the solution {\it without} the need to characterize the spectral functions $A(k)$ and $B(k)$ in terms of the given initial and boundary conditions. This can be achieved by employing the Deift-Zhou method \cite{DZ1992, DZ1993} for the long-time asymptotics \cite{FI1992, FI1992b, FI1994, FI1996} and the Deift-Zhou-Venakides method \cite{DVZ1994, DVZ1997} for the zero-dispersion limit \cite{FK2004, K2003}. However, the complete solution of non-linearizable boundary value problems requires the characterization of $\{A(k), B(k)\}$ in terms of the given initial and boundary conditions. The effective solution of this problem is particularly important for the physically significant case of boundary conditions which are periodic in $t$, since in this case it is {\it not} possible to obtain the rigorous form of the long time asymptotics, without the full characterization of $A(k)$ and $B(k)$ (in spite of this difficulty important results about the large $t$-asymptotics of such problems are presented in \cite{BIK2007, BIK2009, BKSZ2010}). 

In this paper, we revisit the implementation of the unified method to integrable evolution PDEs on the half-line. We concentrate on the {\it effective} characterization of the spectral functions for non-linearizable boundary conditions. We call a characterization of $\{A(k), B(k)\}$ effective if it fulfills the following requirements: $(a)$ In the linear limit, it yields an effective solution of the linearized boundary value problem, i.e. it yields a solution in the form of an integral which involves the transforms of the given initial and boundary conditions. $(b)$ For `small' boundary conditions, it yields an effective perturbative scheme, i.e. it yields an expression in which each term can be computed uniquely in a well-defined recursive scheme.

The effective characterization of $\{A(k), B(k)\}$ presented here is based on an effective characterization of the unknown boundary values in terms of the given initial and boundary conditions, i.e. on an effective characterization of the so-called {\it generalized Dirichlet to Neumann map}. An effective characterization of this map for the Dirichlet problem of the NLS was presented in \cite{BFS2003} by employing the so-called Gelfand-Levitan-Marchenko (GLM) representations. Following this important development, analogous results for the Dirichlet problem of the sine-Gordon, mKdV, and KdV equations were presented in \cite{F2005, TF2008}. The final formulas of \cite{F2005} did {\it not} involve the GLM representations, but their derivation was based on these representations. Here, we rederive the formulas of \cite{F2005} for the Dirichlet problem of the NLS and of the mKdV directly, i.e. {\it without} using the GLM representation, and also present the analogous formulas for the Neumann problem of the NLS as well as for the mKdV equation in the case that either $q_x(0,t)$ or $q_{xx}(0,t)$ are prescribed as boundary conditions.

Our approach uses three ingredients: $(a)$ The large $k$ asymptotics of the eigenfunction $\Phi(t,k)$ defining $\{A(k), B(k)\}$. $(b)$ The global relation {\it and} the equations obtained from the global relation under the transformations which leave invariant the dispersion relation of the associated linearized equation. $(c)$ A perturbative scheme to establish effectiveness. We emphasize that, in general, {\it all} three ingredients are needed, except for the particular case of the NLS, where as explained in appendix \ref{DMSapp}, it is {\it not} necessary to use the invariance of the global relation. In particular, ingredient $(a)$ yields {\it several} possible formulas for the unknown boundary values, so it is absolutely necessary to employ ingredient $(c)$ in order to choose the one that yields an effective solution.

The important idea that the asymptotics of the associated eigenfunctions can be used for the characterization of the unknown boundary values for the NLS was first introduced in \cite{DMS2001} (the authors of \cite{DMS2001} employed the asymptotics of the eigenfunctionÊ $F(t,k)$ instead of $\Phi(t,k)$, but these eigenfunctions are related by equation (\ref{PhiFrelation})).
In other words, the ingredient $(a)$ above was first introduced in \cite{DMS2001}. However, the ingredients $(b)$ and $(c)$ have not been introduced before. Hence, the approach of \cite{DMS2001} {\it cannot} be generalized to other integral nonlinear PDEs, such as the mKdV. Furthermore, since the ingredient $(c)$ was not introduced in \cite{DMS2001}, it is not a priori clear if any of the formulas presented in \cite{DMS2001} yield an effective characterization. It turns out, as explained in appendix \ref{DMSapp}, that the relevant formulas for the NLS {\it are} effective in the case of the half-line but {\it not} for the case of the finite interval \cite{LFinterval}.

From the above discussion it follows that it is possible to construct an effective characterization of both the generalized Dirichlet to Neumann map and of $\{A(k), B(k)\}$, {\it without} the need to introduce the GLM representations. However, it appears that the latter representations might have an advantage for analyzing the large $t$-asymptotics of the solution in the case of $t$-periodic boundary conditions. For this reason, we revisit these representations in \cite{LFtperiodic}, where we present a simplification, as well as a significant extension, of the results of \cite{BFS2003, F2005}.

\subsection{Linear evolution PDEs with a second versus a third order $x$-derivative}
Our definition of an effective characterization of $A(k)$ and $B(k)$, suggests that before attempting to solve a nonlinear boundary value problem, it is imperative that we have first constructed an effective solution of the underlying linearized problem. In this respect we note that linear evolution PDEs with second order derivatives are rather special, namely they can be solved by certain transforms which {\it cannot} be applied to linear evolution PDEs with third or higher order derivatives. In this sense, there exists a crucial difference between the linearized version of the NLS, namely of the equation
\begin{equation}\label{linearizedNLS}
  iu_t + u_{xx} = 0,
\end{equation}
and the linearized version of the mKdV, namely of the equation
\begin{equation}\label{linearizedmKdVI}
  u_t + u_{xxx} = 0.
\end{equation}
Let us consider the Dirichlet problem on the half-line for equations (\ref{linearizedNLS}) and (\ref{linearizedmKdVI}), i.e. let us consider (\ref{linearizedNLS}) and (\ref{linearizedmKdVI}) in the domain
\begin{equation}\label{halflinedomain}
   0 < x < \infty; \qquad 0 < t < T, \qquad T \text{ positive constant},
\end{equation}
with
\begin{equation}\label{u0g0data}
  u(x,0) = u_0(x), \quad 0 < x < \infty; \qquad u(0,t) = g_0(t), \quad 0< t < T,
\end{equation}
where $u_0(x)$ and $g_0(t)$ are given functions with appropriate smoothness and decay, satisfying $u_0(0) = g_0(0)$. It is well known that the classical sine transform can be used for solving equation (\ref{linearizedNLS}) with the conditions (\ref{u0g0data}).
We recall that this transform can be obtained from the application of the Fourier transform to an odd function. Thus, the application of the sine transform is {\it equivalent} to solving the above problem via an odd extension of $u$ from the half-line to the full line. Similarly, it is possible to solve equation (\ref{linearizedNLS}) on the half-line with Neumann boundary conditions, either by using the classical cosine-transform, or equivalently via an even extension of $u$ from the half-line to the full line. The solution of the Robin boundary value problem can also be constructed by either using an appropriate $x$-transform, or via a certain extension of $u$ from the half-line to the full line. 

In contrast to equation (\ref{linearizedNLS}), there does {\it not} exist a classical $x$-transform for solving equation (\ref{linearizedmKdVI}) on the half-line. 

The unified method provides an elegant and effective solution for {\it both} equations (\ref{linearizedNLS}) and (\ref{linearizedmKdVI}). Indeed, define the following transforms of the given functions $u_0(x)$ and $g_0(t)$,
$$\hat{u}_0(k) = \int_0^\infty e^{-ikx} u_0(x) dx, \quad \text{Im}\;k \leq 0; \qquad
\tilde{g}_0(k) = \int_0^T e^{k\tau} g_0(\tau) d\tau, \quad k \in \C.
$$
Then,
\begin{equation}\label{urepresentation}
u(x,t) = \frac{1}{2\pi}\int_{-\infty}^\infty e^{ikx - i\omega(k)t}\hat{u}_0(k) dk - \frac{1}{2\pi}\int_{\partial D^+} e^{ikx - i\omega(k)t} \tilde{g}(k) dk,
\end{equation}
where the $(x,t)$-domain is defined in (\ref{halflinedomain}), $\omega(k)$ is the dispersion relation, $\partial D^+$ is the oriented boundary of the domain $D^+$, where $D^+$ is defined by
$$D^+ = \{\text{Im}\;k \geq 0, \; \text{Im}\; \omega(k) > 0\},$$
with the assumption that $D^+$ is to the left of the increasing orientation, see Figure \ref{partialDplusfig}, and $\tilde{g}(k)$ can be explicitly computed in terms of $\hat{u}_0$ and $\tilde{g}_0$:
\begin{figure}
\begin{center}
\bigskip \bigskip
 \begin{overpic}[width=.4\textwidth]{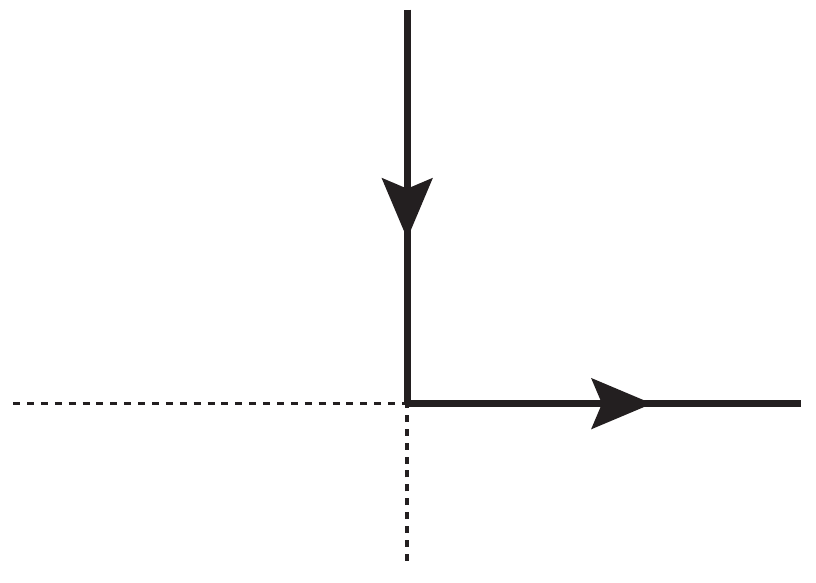}
      \put(0,70){$(a)$}
      \put(100,20){$\text{Re}\; k$}
      \put(43,73){$\text{Im}\; k$}
    \end{overpic}
    \qquad \qquad
\begin{overpic}[width=.4\textwidth]{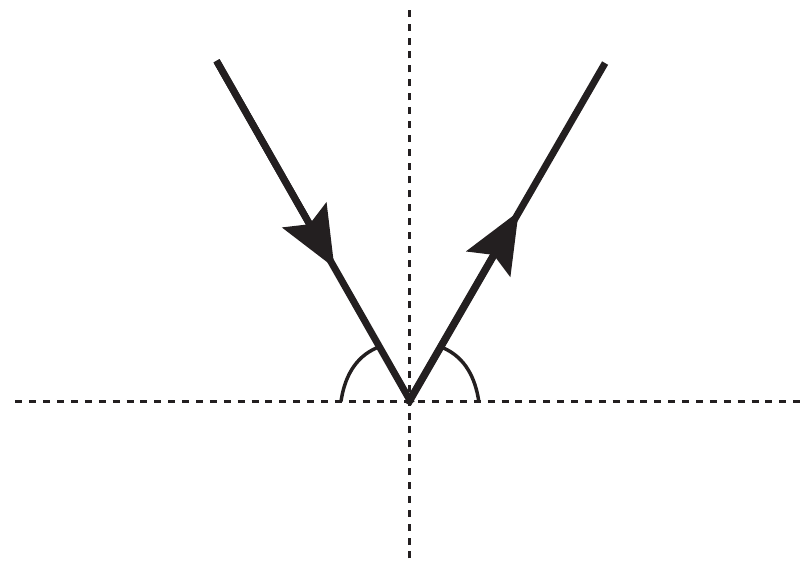}
      \put(0,70){$(b)$}
       \put(101,19){$\text{Re}\; k$}
      \put(43,73){$\text{Im}\; k$}
      \put(60,25){$\pi/3$}
      \put(30,25){$\pi/3$}
      \end{overpic}
     \begin{figuretext}\label{partialDplusfig}
       The contour $\partial D^+$ in the complex $k$-plane for equation (\ref{linearizedNLS}), figure $(a)$, and equation (\ref{linearizedmKdVI}), figure $(b)$.
     \end{figuretext}
     \end{center}
\end{figure}

\medskip
\noindent
\underline{Linearized NLS}
\begin{equation}\label{nlsdispersion} 
 \omega(k) = k^2, \qquad \tilde{g}(k) = \hat{u}_0(-k) - 2k\tilde{g}_0(i k^2), \qquad \text{Im}\; k \geq 0.
\end{equation}

\noindent
\underline{Linearized mKdV}
$$\omega(k) = -k^3, \qquad \tilde{g}(k) = - \alpha \hat{u}_0(\alpha k) - \alpha^2 \hat{u}_0(\alpha^2k) + 3k^2\tilde{g}_0(-ik^3), \qquad \alpha = e^{\frac{2i\pi}{3}}, \; \text{Im}\; k \geq 0.$$

By deforming the contour of Figure \ref{partialDplusfig}$a$ from the positive imaginary axis to the negative real axis, it is possible to recover from equation (\ref{urepresentation}) the classical solution obtained via the sine-transform (or equivalently via an odd extension):
\begin{equation}\label{urepclassical}
u(x,t) = \frac{2}{\pi}\int_{0}^\infty  e^{-ik^2t} \sin{kx} \left[\int_0^\infty \sin(k\xi) u_0(\xi) d\xi + ik\int_0^t e^{ik^2\tau}g_0(\tau) d\tau\right]dk.
\end{equation}
However, as noted earlier, the right-hand side of equation (\ref{urepclassical}) is {\it not} uniformly convergent at $x = 0$. Indeed, if it were uniformly convergent then one could take the limit $x \to 0$ inside the integral and then $u(0,t)$ would vanish, whereas we know that $u(0,t) = g_0(t)$. {\it This lack of uniform convergence is a significant disadvantage of the solution obtained via a classical transform method for any boundary value problem.} This fact, in addition to rendering such representations unsuitable for numerical computations, also makes it difficult to prove existence without using PDE techniques. Indeed, taking into consideration that the construction of the solution via any transform method {\it assumes} existence of solutions, one must verify that the solution constructed under this a priori assumption satisfies the given PDE and the given initial and boundary conditions. An important advantage of the new method is that in addition to yielding an effective representation for evolution PDEs with spatial derivatives of arbitrary order, it also yields integral representations which are uniformly convergent on the boundary. This makes it straightforward to verify that the integral representation indeed solves the given initial-boundary value problem. For example, substituting (\ref{nlsdispersion}) into (\ref{urepresentation}) and evaluating the resulting equation at $x = 0$, we find
\begin{align}\nonumber
u(0,t) = &\; \frac{1}{2\pi} \int_{-\infty}^\infty e^{-ik^2t}\hat{u}_0(k)dk - \frac{1}{2\pi} \int_{\partial D^+}e^{-ik^2t}\hat{u}_0(-k) dk
	\\ \label{linearnlsrep}
& + \frac{1}{2\pi} \int_{\partial D^+} 2ke^{-ik^2t}\left(\int_0^Te^{ik^2\tau}g_0(\tau)d\tau\right)dk.
\end{align}
By deforming the contour of integration in the second integral on the right-hand side of (\ref{linearnlsrep}) from the positive imaginary axis to the negative real axis and then replacing $k$ with $-k$ in the resulting integral, it follows that the first two integrals on the right-hand side of (\ref{linearnlsrep}) cancel; by making the change of variables $k^2 = l$ in the third integral on the right-hand side of (\ref{linearnlsrep}) and then using the integral formula of the classical Fourier transform, it follows that $u(0,t) = g_0(t)$.

\subsection{Organization of the paper}
In section \ref{prelsec} we review the main steps of the unified method; proofs of the basic results can be found in \cite{F2002, FIS2005}. In sections \ref{NLSsec} and \ref{mKdVIsec} we derive effective characterizations of the spectral functions $A(k)$ and $B(k)$ for the Dirichlet and Neumann boundary value problems for the NLS and mKdV equations, respectively. In section \ref{D2Nsec} we utilize these characterizations to give a perturbative construction of the Dirichlet to Neumann map. In section \ref{conclusionsec} the above results are discussed further.

\subsection{Basic assumptions and notations}
\begin{itemize}
\item We assume that $q(x,t)$ vanishes sufficiently fast for all $t$ as $x \to \infty$. Furthermore, we assume that the given initial conditions satisfy
\begin{align}
  q_0(x) := q(x,0) \in L_1(\R^+) \cap L_2(\R^+).
\end{align}

\item $g_0$, $g_1$, $g_2$ will denote $q$, $q_x$, and $q_{xx}$ evaluated at $x = 0$, i.e.
\begin{align}
  g_0(t) = q(0,t), \qquad g_1(t) = q_x(0,t), \qquad g_2(t) = q_{xx}(0,t), \qquad 0 < t < T.
\end{align}

\item Let $A_1$ and $A_2$ denote the two column vectors of the $2 \times 2$ matrix $A$. The notation $A(k)$, $k \in (D_1, D_2)$, means that for $A_1(k)$, $k \in D_1$, and for $A_2(k)$, $k \in D_2$.
\end{itemize}

\section{Preliminaries}\label{prelsec}\nequation
Several important integrable nonlinear PDEs admit the following Lax pair \cite{Lax1968} formulation:
\begin{subequations}\label{generallax}
\begin{align}\label{generallaxa}
 & \frac{\partial \mu}{\partial x}(x,t,k) + if_1(k) \hat{\sigma}_3 \mu(x,t,k) = Q(x,t,k) \mu(x,t,k),
  	\\\label{generallaxb}
&  \frac{\partial \mu}{\partial t}(x,t,k) + if_2(k) \hat{\sigma}_3 \mu(x,t,k) = \tilde{Q}(x,t,k) \mu(x,t,k), \qquad x,t \in \R, \quad k \in \C,
\end{align}
\end{subequations}
where $\mu$ is a $2 \times 2$ matrix-valued function, $\{f_1(k), f_2(k)\}$ are given analytic functions of $k$, and the $2 \times 2$ matrix-valued functions $Q$ and $\tilde{Q}$ are given analytic function of $k$, of $q(x,t)$, of $\bar{q}(x,t)$, and of the $x$-derivatives of these functions. The action of $\hat{\sigma}_3$ on a $2 \times 2$-matrix $A$ is defined by
\begin{subequations}
\begin{align}
  \hat{\sigma}_3A = [\sigma_3, A], \qquad \sigma_3 = \text{diag}(1,-1),
\end{align}
and hence
\begin{align}\label{ehatsigma3}
  e^{\hat{\sigma}_3x} A = e^{\sigma_3x}Ae^{-\sigma_3x} = \begin{pmatrix} A_{11} & e^{2x}A_{12} \\
  e^{-2x} A_{21} & A_{22} \end{pmatrix}.
\end{align}
\end{subequations}

\medskip
\noindent
\underline{Example 1}
For the NLS equation,
\begin{subequations}\label{NLSQQtildedef}
\begin{align}\label{nls}
 & i\frac{\partial q}{\partial t} + \frac{\partial^2 q}{\partial x^2} - 2 \lambda |q|^2 q = 0, \qquad \lambda = \pm 1,
  	\\
&  f_1(k) = k, \qquad f_2(k) = 2k^2, \qquad Q(x,t) = \begin{pmatrix} 0 & q(x,t) \\ \lambda \bar{q}(x,t) & 0 \end{pmatrix},
  	\\
&  \tilde{Q}(x,t) = 2kQ(x,t) + \tilde{Q}^{(1)}(x,t), \qquad \tilde{Q}^{(1)}(x,t) = -i(Q_x(x,t) + \lambda |q(x,t)|^2 )\sigma_3.
\end{align}
\end{subequations}

\medskip
\noindent
\underline{Example 2}
For the mKdV equation,
\begin{subequations}
\begin{align}\label{mkdvI}
 & \frac{\partial q}{\partial t} + \frac{\partial^3 q}{\partial x^3} - 6 \lambda q^2 \frac{\partial q}{\partial x} = 0, \qquad q\;\text{real},\quad \lambda = \pm 1,
  	\\
&  f_1(k) = k, \qquad f_2(k) = 4k^3, \qquad Q(x,t) = \begin{pmatrix} 0 & q(x,t) \\ \lambda q(x,t) & 0 \end{pmatrix},
  	\\
&  \tilde{Q}(x,t) = 2Q(x,t)^3 - Q_{xx}(x,t) - 2ik(Q(x,t)^2 + Q_x(x,t))\sigma_3 + 4k^2Q(x,t).
\end{align}
\end{subequations}

\subsection{The basic eigenfunctions}
Using the identity (\ref{ehatsigma3}), it follows that equations (\ref{generallax}) can be written in the form 
\begin{subequations}
\begin{align}
&  \frac{\partial}{\partial x} e^{i(f_1(k)x + f_2(k)t)\hat{\sigma}_3}\mu(x,t,k) = e^{i(f_1(k)x + f_2(k)t)\hat{\sigma}_3}Q(x,t,k)\mu(x,t,k),
  	\\	
&  \frac{\partial}{\partial t} e^{i(f_1(k)x + f_2(k)t)\hat{\sigma}_3}\mu(x,t,k) = e^{i(f_1(k)x + f_2(k)t)\hat{\sigma}_3}\tilde{Q}(x,t,k)\mu(x,t,k), 
	\\ \nonumber
& \hspace{8cm} x,t \in \R, \quad k \in \C.
\end{align}
\end{subequations}
Hence,
\begin{subequations}
\begin{align}\label{generallaxdiffform}
  d\Bigl[e^{i(f_1(k)x + f_2(k)t)\hat{\sigma}_3}\mu(x,t,k)\Bigr] = e^{i(f_1(k)x + f_2(k)t)\hat{\sigma}_3}W(x,t,k), \qquad x,t \in \R, \; k \in \C,
\end{align}
where the differential form $W$ is defined by
\begin{align}\label{Wdef}
  W(x,t,k) = \left(Q(x,t,k) dx + \tilde{Q}(x,t,k)dt\right) \mu(x,t,k).
\end{align}
\end{subequations}
Equation (\ref{generallaxdiffform}) implies that if $(x,t) \in \Omega$, where $\Omega \subset \R^2$ is a simply connected domain, then the following equation characterizes a function $\mu_j(x,t,k)$ which satisfies {\it both} (\ref{generallaxa}) and (\ref{generallaxb}):
\begin{align}\label{mujdef}
  \mu_j(x,t,k) = I + \int_{(x_j, t_j)}^{(x,t)} e^{-i(f_1(k)(x - \xi)\hat{\sigma}_3 - if_2(k)(t - \tau)\hat{\sigma}_3} W_j(\xi,\tau,k), \quad
  (x,t), (x_j, t_j) \in \Omega,
\end{align}
where $W_j$ is the differential form defined in (\ref{Wdef}) with $\mu$ replaced by $\mu_j$ and $(x_j, t_j)$ is a fixed point in $\Omega$. 

The implementation of the unified method relies on the fact that if $\Omega$ is a convex polygon, then by choosing $(x_j, t_j)$ as the collection of the vertices of this polygon, it follows that the collection of $\mu_j$ defines a sectionally analytic function in the complex $k$-plane \cite{F1997}.

In the case of the half-line, $\Omega$ is given by
\begin{align}\label{Omegadef}  
  \Omega = \{0 < x < \infty, \; 0 < t < T\},
\end{align}
thus the relevant polygon is convex and the three corners are the points
$$(0,T), \quad (0,0), \quad (\infty,t).$$
\begin{figure}
\begin{center}
\medskip
 \begin{overpic}[width=.5\textwidth]{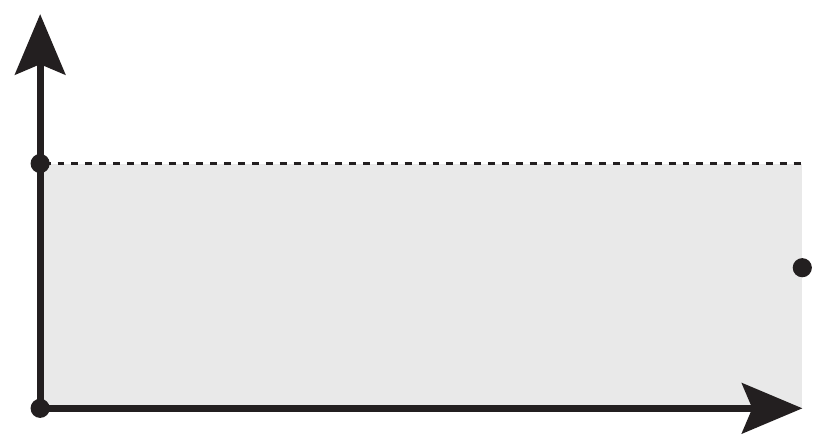}
      \put(8,36){$1$}
      \put(-1,0){$2$}
      \put(101,20){$3$}
      \put(100,3){$x$}
      \put(-2.5,32.5){$T$}
      \put(48,18){$\Omega$}
    \end{overpic}
    \qquad \qquad
     \begin{figuretext}\label{halflinedomain.pdf}
       The domain $\Omega$ with the corners at $(0,T)$, $(0,0)$, and $(\infty,t)$.
     \end{figuretext}
     \end{center}
\end{figure}
We refer to these corners as $\{1,2,3\}$, respectively.
\begin{figure}
\begin{center}
\medskip
 \begin{overpic}[width=.25\textwidth]{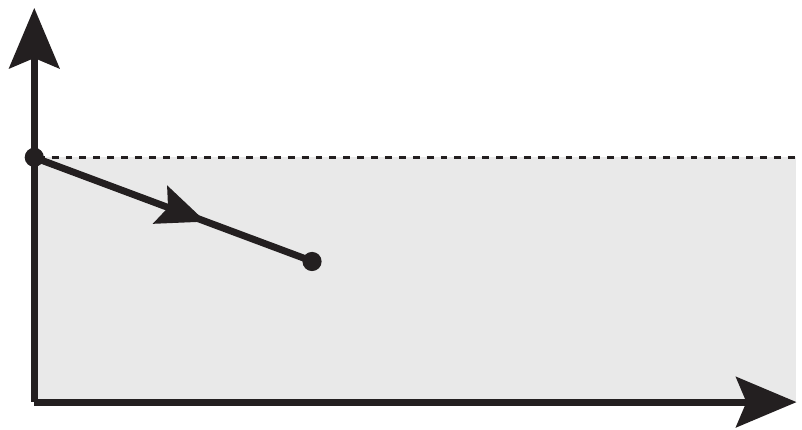}
 \put(-8,30){\small $T$}
 \put(43,17){\small $(x,t)$}
 \put(48,-10){\small (1)}
   \end{overpic}
   \qquad\quad
  \begin{overpic}[width=.25\textwidth]{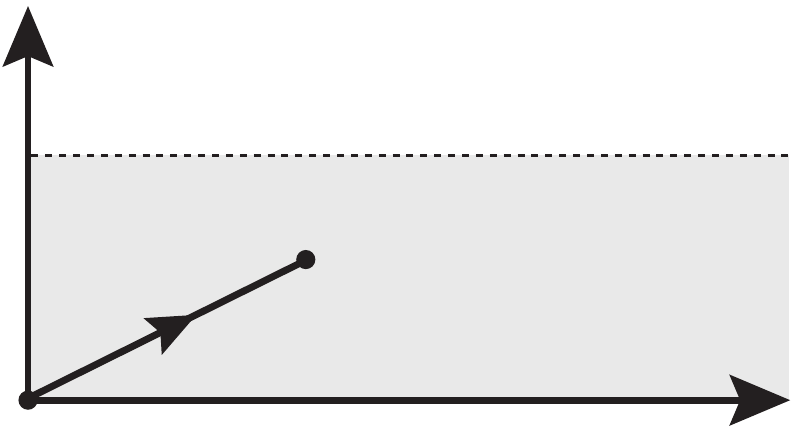}
    \put(-8,30){\small $T$}
    \put(43,18){\small $(x,t)$} 
 \put(48,-10){\small (2)}
  \end{overpic} 
  \qquad \quad
  \begin{overpic}[width=.25\textwidth]{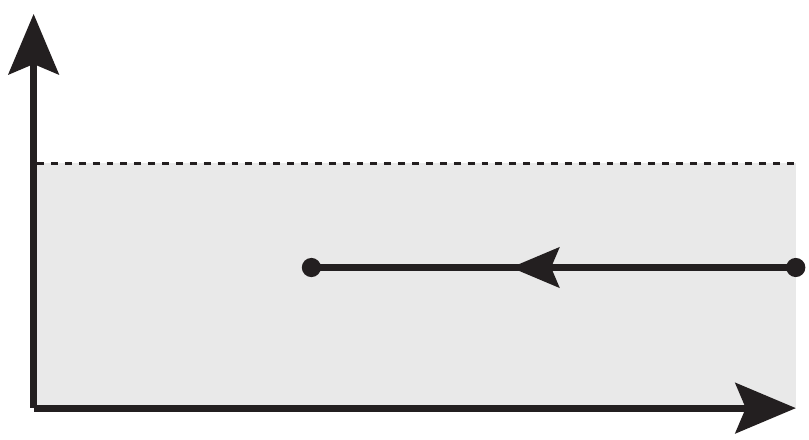}
     \put(-8,30){\small $T$}
    \put(13,18){\small $(x,t)$} 
 \put(48,-10){\small (3)}
  \end{overpic}
  \bigskip
   \begin{figuretext}\label{mucontoursfig}
     The contours used to define $\{\mu_j\}_1^3$.
   \end{figuretext}
   \end{center}
\end{figure}
Let $\{\mu_1, \mu_2, \mu_3\}$ be the eigenfunctions corresponding to these corners. These eigenfunctions are {\it independent} of the contours from $(x_j, t_j)$ to $(x,t)$. Typical contours are shown in figure \ref{mucontoursfig}. 

For $\mu_1, \mu_2$, $\mu_3$, the following inequalities are valid respectively:
\begin{align*}
   \begin{cases} x - \xi > 0 \\ t - \tau < 0 \end{cases},  \qquad 
   \begin{cases} x - \xi > 0 \\ t - \tau > 0 \end{cases},  \qquad 
    \begin{cases} x - \xi < 0 \\ t - \tau = 0 \end{cases}.
\end{align*}
Hence, the expressions that appear in the first and second columns of (\ref{mujdef}) are analytic and bounded in the following domains of the complex $k$-plane: 
\begin{align*}
 & \mu_1: \; \left(\text{Im}\; f_1 > 0 \cap \text{Im}\; f_2 < 0, \; \text{Im}\; f_1 < 0 \cap \text{Im}\; f_2 > 0\right) =: (D_2, D_3),
  	\\
 & \mu_2: \; \left(\text{Im}\; f_1 > 0 \cap \text{Im}\; f_2 > 0, \; \text{Im}\; f_1 < 0 \cap \text{Im}\; f_2 < 0\right) =: (D_1, D_4),
	\\
&  \mu_3: \; \left(\text{Im}\; f_1 < 0,\; \text{Im}\; f_1 > 0\right) =: (\C^-, \C^+).
\end{align*}
\begin{figure}
\begin{center}
\bigskip \bigskip
 \begin{overpic}[width=.29\textwidth]{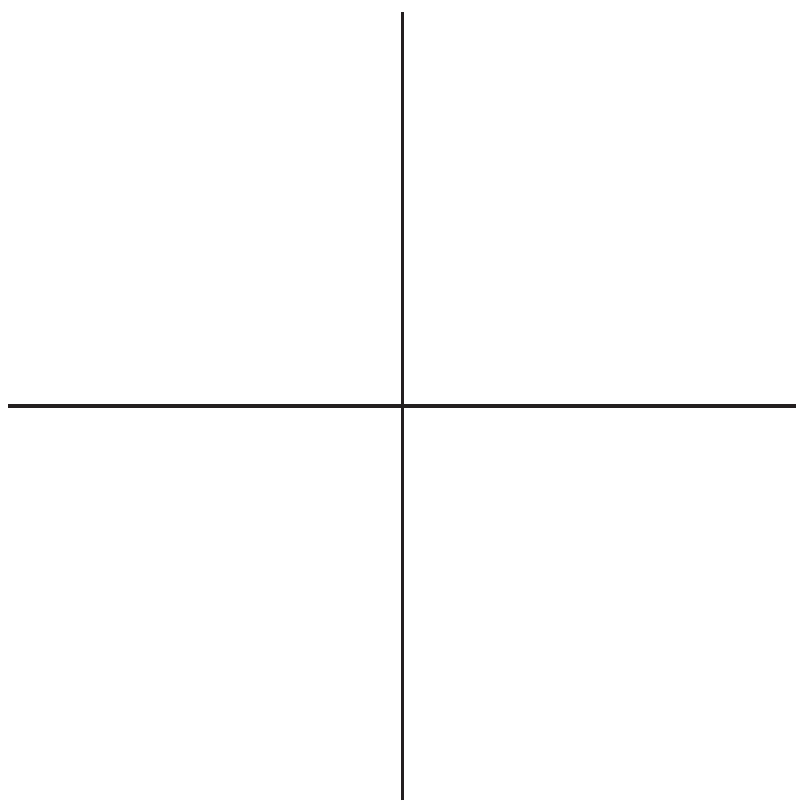}
      \put(-4,94){$(a)$}
      \put(100,47){$\text{Re}\; k$}
      \put(68,70){$D_1$}
      \put(20,70){$D_2$}
       \put(68,25){$D_4$}
      \put(20,25){$D_3$}
    \end{overpic}
    \qquad \qquad
\begin{overpic}[width=.34\textwidth]{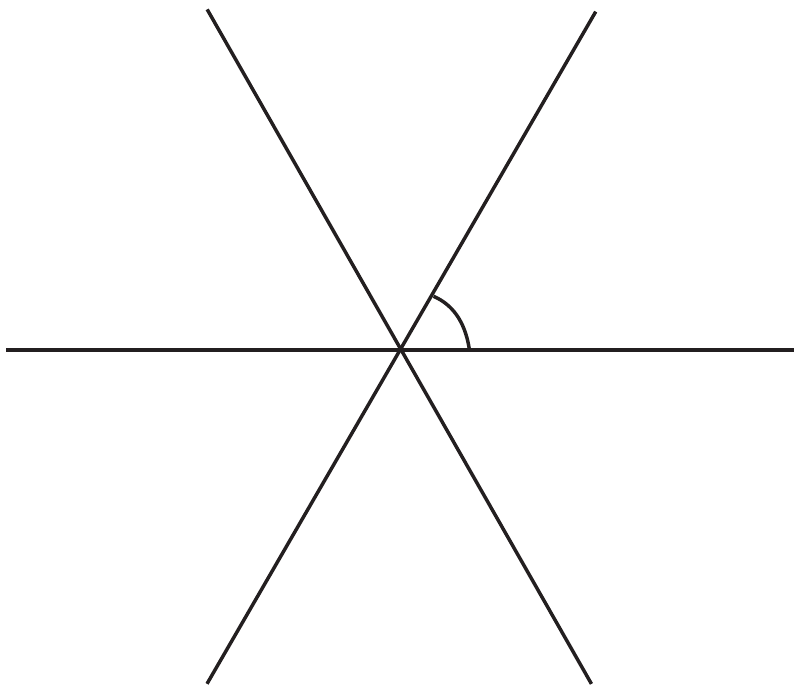}
      \put(-5,80){$(b)$}
      \put(77,57){$D_1$}
      \put(46,70){$D_2$}
      \put(15,57){$D_1$}
      \put(15,25){$D_4$}
      \put(46,12){$D_3$}
      \put(77,25){$D_4$}
       \put(101,40){$\text{Re}\; k$}
      \put(60,47){$\pi/3$}
      \end{overpic}
     \begin{figuretext}\label{Djsfig}
       The domains $\{D_j\}_1^4$ for $(a)$ NLS and $(b)$ mKdV.
     \end{figuretext}
     \end{center}
\end{figure}
For NLS and mKdV the domains $\{D_j\}_1^4$ are shown in figure \ref{Djsfig}. For these equations, $f_1(k) = k$ so that $\C^-$ and $\C^+$ denote the lower and upper halves of the complex $k$-plane, respectively.

\subsection{The spectral functions}
The crucial advantage of requiring that the functions $\{\mu_j\}_1^3$ solve equation (\ref{mujdef}) is that any two such functions are related via a matrix which has {\it explicit exponential} $(x,t)$-dependence of the form $\exp[-if_1(k)x\hat{\sigma}_3 - if_2(k)t\hat{\sigma}_3]\rho(k)$, where the function $\rho(k)$ can be computed by evaluating the relevant relation at any convenient point in the domain $\Omega$. In particular, using that $\mu_2(0,0,k) = I$, we find the following relations:
\begin{subequations}\label{mu321relations}
\begin{align}
  \mu_3(x,t,k) = \mu_2(x,t,k)e^{-i(f_1(k)x + f_2(k) t)\hat{\sigma}_3}s(k), \qquad k \in (\C^-, \C^+),
  	\\
  \mu_1(x,t,k) = \mu_2(x,t,k)e^{-i(f_1(k)x + f_2(k) t)\hat{\sigma}_3}S(k), \qquad k \in (D_2, D_3),	
\end{align} 
\end{subequations}
where the spectral functions $s(k)$ and $S(k)$ are defined by
\begin{align}
 & s(k) = \Psi(0,k), \qquad k \in (\C^-, \C^+),
  	\\
&  S^{-1}(k) = e^{if_2T \hat{\sigma}_3}\Phi(T,k), \qquad k \in \C,
\end{align}
and $\Psi$, $\Phi$ denote the functions obtained by evaluating $\mu_3$ and $\mu_2$ at $t = 0$ and $x = 0$ respectively, i.e.
$$\Psi(x,k) = \mu_3(x,0,k), \qquad k \in (\C^- , \C^+); \qquad \Phi(t,k) = \mu_2(0,t,k), \qquad k \in \C.$$
Hence,Ê $\Psi$ and $\Phi$ satisfy the $x$ and $t$-parts of the associated Lax pair evaluated at $t = 0$ and $x = 0$ respectively:
\begin{align}\nonumber
&  \Psi_x + if_1 \hat{\sigma}_3 \Psi = Q_0 \Psi, \qquad Q_0(x,k) = Q(x,0,k), 
  	\\ \label{Psieq}
&  \lim_{x \to \infty} \Psi = I, \qquad 0 < x < \infty, \quad k \in (\C^-, \C^+),
\end{align}
and
\begin{align} \nonumber
&  \Phi_t + if_2 \hat{\sigma}_3 \Phi = \tilde{Q}_0 \Phi, \qquad \tilde{Q}_0(t,k) = \tilde{Q}(0,t,k),
  	\\ \label{Phieq}
&  \Phi(0,t) = I, \qquad 0 < t < T, \quad k \in \C.
\end{align}
Equations (\ref{Psieq}) and (\ref{Phieq}) are equivalent to the following linear Volterra integral equations:
\begin{align}\label{Psivolterra}
  \Psi(x,k) = I - \int_x^\infty e^{-if_1(k)(x - \xi)\hat{\sigma}_3} (Q_0\Psi)(\xi, k) d\xi, \qquad 0 < x < \infty, \quad k \in (\C^-, \C^+)
\end{align}
and
\begin{align}\label{Phivolterra}
  \Phi(t,k) = I + \int_0^t e^{-if_2(k)(t - \tau)\hat{\sigma}_3} (\tilde{Q}_0\Phi)(\tau, k) d\tau, \qquad 0 < t < T, \quad k \in \C.
\end{align}

\subsection{The global relation}
Let $F(t,k)$ denote the eigenfunction obtained from the evaluation of $\mu_3$ at $x = 0$, i.e.
$$F(t,k) = \mu_3(0,t,k), \qquad k \in (\C^-, \C^+).$$
The eigenfunctions $\Phi$ and $F$ satisfy the same differential equation (\ref{generallaxb}), thus they are simply related:
\begin{align}\label{PhiFrelation}
  \Phi(t,k)e^{-if_2(k) t\hat{\sigma}_3}s(k) = F(t,k), \qquad 0 < t < T, \quad k \in (\C^-, \C^+).
\end{align}
Evaluating this equation at $t = T$ we obtain the {\it global relation}
\begin{align}  \label{globalrelation}
  S^{-1}(k)s(k) = e^{if_2T \hat{\sigma}_3}F(T,k), \qquad k \in (\C^-, \C^+).
\end{align}

\subsection{Notations}
The superscripts $(+,-,1,2,3,4)$ will denote the domain of analyticity of the relevant vector eigenfunction:
\begin{align}
  \mu_1 = \left(\mu_1^{(2)}, \mu_1^{(3)}\right), \qquad
   \mu_2 = \left(\mu_2^{(1)}, \mu_2^{(4)}\right), \qquad
    \mu_3 = \left(\mu_3^{-}, \mu_3^{+}\right).
\end{align}
For NLS and mKdV, the matrices $Q$ and $\tilde{Q}$ possess certain symmetries, which in turn imply certain symmetries for the eigenfunctions and for the spectral functions:
\begin{align}\label{ABdef}
 & \Phi(t,k) = \begin{pmatrix}\overline{\Phi_2(t,\bar{k})} & \Phi_1(t,k) \\
  \lambda \overline{\Phi_1(t, \bar{k})} & \Phi_2(t,k) \end{pmatrix}, \qquad 
  S(k) = \begin{pmatrix}\overline{A(\bar{k})} & B(k) \\
  \lambda \overline{B(\bar{k})} & A(k) \end{pmatrix}, \qquad k \in \C,
	\\ \label{abdef}
&  \Psi(x,k) = \begin{pmatrix}\overline{\Psi_2(x,\bar{k})} & \Psi_1(x,k) \\
  \lambda \overline{\Psi_1(x, \bar{k})} & \Psi_2(x,k) \end{pmatrix}, \qquad 
  s(k) = \begin{pmatrix}\overline{a(\bar{k})} & b(k) \\
  \lambda \overline{b(\bar{k})} & a(k) \end{pmatrix}, \qquad k \in (\C^-, \C^+).
\end{align}
We will henceforth assume that these symmetries are present.

\subsection{The RH problem}\label{RHsubsec}
Equations (\ref{mu321relations}) together with the estimate
\begin{align}
  \mu_j = I + O\left(\frac{1}{k}\right), \qquad k \to \infty, \quad j = 1,2,3,
\end{align}
for $k$ in the domain in which $\mu_j$ is bounded and analytic, can be used to formulate a RH problem. Indeed, equations (\ref{mu321relations}) can be written in the form
\begin{align}\label{jumpcondition}
  M_-(x,t,k) = M_+(x,t,k) J(x,t,k), \qquad k \in \bar{D}_+ \cap \bar{D}_-,
\end{align}
where
$$D_+ = D_1 \cup D_3, \qquad D_- = D_2 \cup D_4,$$ 
and the matrices $M_-$, $M_+$, $J$ are defined as follows:
\begin{align*}
& M_+ = \left(\frac{\mu_2^{(1)}}{a(k)}, \mu_3^+\right), \quad k \in D_1; \qquad
M_- = \left(\frac{\mu_1^{(2)}}{d(k)}, \mu_3^{+}\right), \quad k \in D_2; \qquad
	\\
& M_+ = \left(\mu_3^{-}, \frac{\mu_1^{(3)}}{\overline{d(\bar{k})}} \right), \quad k \in D_3; \qquad
M_- = \left(\mu_3^{-}, \frac{\mu_2^{(4)}}{\overline{a(\bar{k})}} \right), \quad k \in D_4; 
	\\
&J(x,t,k) = \begin{cases} J_1, & k \in \bar{D}_1 \cap \bar{D}_2,	\\
J_2 = J_3J_4^{-1}J_1, & k \in \bar{D}_2 \cap \bar{D}_3,	\\
J_3, & k \in \bar{D}_3 \cap \bar{D}_4,	\\
J_4, & k \in \bar{D}_4 \cap \bar{D}_1;
\end{cases}
	\\
& J_1 = \begin{pmatrix} 1 & 0 \\ \frac{\lambda \overline{B(\bar{k})}}{a(k) d(k)} e^{2i\theta} & 1 \end{pmatrix}, \quad
J_3 =  \begin{pmatrix} 1 & -\frac{B(k)}{\overline{a(\bar{k})} \overline{d(\bar{k})}} e^{-2i\theta} \\ 0 & 1 \end{pmatrix}, \quad
J_4 = \begin{pmatrix} 1 & -\frac{b(k)}{\overline{a(\bar{k})}}e^{-2i\theta} \\
\frac{\lambda \overline{b(\bar{k})}}{a(k)}e^{2i\theta} & \frac{1}{a(k) \overline{a(\bar{k})}} \end{pmatrix},
\end{align*}
with
$$d(k) = a(k)\overline{A(\bar{k})} - \lambda b(k) \overline{B(\bar{k})}; \qquad \theta(x,t,k) = f_1(k)x + f_2(k)t.$$

\subsection{The main rigorous result}
The above derivation is based on the assumption that there exists a solution of the given nonlinear PDE with appropriate smoothness and decay. The rigorous justification of the unified transform method is based on the following important result, first proved in \cite{FIS2005} in connection with the NLS equation:
(i) Given initial condition $q_0(x) \in L_1(\R^+) \cap L_2(\R^+)$, define the spectral functions $\{a(k), b(k)\}$ by 
$$a(k) = \Psi_2(0,k), \qquad b(k) = \Psi_1(0,k), \qquad k \in \C^+,$$
where $\Psi(x,k)$ is the solution of (\ref{Psieq}), or equivalently of (\ref{Psivolterra}) and $Q_0(x,k)$ is uniquely defined in terms of $q_0(x)$ (with $q(x,0)$ replaced by $q_0(x)$).
(ii) Given a subset of the boundary values 
$$g_0(t) = q(0,t), \qquad g_1(t) = q_x(0,t), \qquad g_2(t) = q_{xx}(0,t),$$
{\it assume} that it is possible to construct the spectral functions $\{A(k), B(k)\}$,
\begin{align}
  A(k) = \overline{\Phi_2(t, \bar{k})}, \qquad B(k) = -e^{2if_2T}\Phi_1(T,k), \qquad k \in \C,
\end{align}
under the following requirements: $(a)$ $\Phi(t,k)$ is a solution of (\ref{Phieq}), or equivalently of (\ref{Phivolterra}), where $\tilde{Q}_0(t,k)$ is defined in terms of $g_j$, $j = 0,1,2$. $(b)$ $\Phi$ satisfies the restriction
\begin{align}\label{PhisF}
  \Phi(t,k)e^{-if_2(k) t\hat{\sigma}_3}s(k) = F(t,k), \qquad 0 < t < T, \quad k \in (\C^-, \C^+),
\end{align}
where $F(t,k)$ is some function whose first and second columns are analytic in k in $\C^-$ and $\C^+$ and
$$F(t,k) = I + O\biggl(\frac{1}{k}\biggr), \qquad k \to \infty, \quad k \in (\C^-, \C^+).$$
(iii) Define the RH problem with the jump condition (\ref{jumpcondition}) in terms of the spectral functions $\{a(k), b(k), A(k), B(k)\}$. Then, the unique solution of this RH problem can be used to compute a function $q(x,t)$ which solves the given nonlinear PDE and furthermore,
\begin{align}
  q(x,0) = q_0(x), \qquad 0 < x < \infty; \quad \frac{\partial^j q}{\partial x^j}(0,t) = g_j(t), \qquad j = 0,1,2, \quad 0 < t < T.
\end{align}

The above theorem shows that the rigorous justification of the unified method reduces to the following problem:
Given $\{a(k), b(k)\}$ and a subset of $\{g_0, g_1, g_2\}$, characterize $\{A(k), B(k)\}$ from the requirement that $\Phi(t,k)$ satisfies (\ref{Phieq}) and (\ref{PhisF}). 

It turns out that:  $(a)$ For linearizable boundary conditions, it is possible to express {\it explicitly} $\{A(k), B(k)\}$ in terms of $\{a(k), b(k)\}$ and the given boundary conditions. $(b)$ For non-linearizable boundary conditions, it is possible to express $\{A(k), B(k)\}$ in terms of a system of two nonlinear equations which can be iteratively solved for `small' boundary conditions.

\begin{remark} \upshape
  For the NLS, there exists {\it one} unknown boundary value. Thus, one might expect that since the global relation (\ref{PhisF}) provides {\it one} equation connecting this unknown function with the given initial and boundary conditions, in this case it is possible, by utilizing the global relation, to characterize the unknown boundary value. However, for the mKdV, there exist {\it two} unknown boundary values. In this case, the solution of the associated linear problem suggests that it is {\it impossible} to solve this problem, unless one uses the transformation which leaves $f_2(k)$ invariant in order to obtain an {\it additional} equation from the global relation.
\end{remark}

\begin{remark}\upshape
  It is important to note that the rigorous result mentioned earlier does {\it not} require the knowledge of the explicit form of $F$; it only requires the \textit{existence} of a function $F$ with specific analyticity properties. This suggests that the most efficient approach of characterizing $\{A(k), B(k)\}$ is to actually {\it eliminate} $F$. This is precisely the philosophy used in the linear limit for the derivation of the generalized Dirichlet to Neumann map and it is also the philosophy used in \cite{BFS2003, F2005, TF2008}.
\end{remark}

\section{The NLS equation}\label{NLSsec}\nequation
We consider the Dirichlet and Neumann boundary value problems for the NLS equation (\ref{nls}) posed on the half-line.

\subsection{The global relation}
Define $c(t, k)$ by 
\begin{equation}\label{cddef}
  c(t, k) = \frac{1}{a(k)}(F(t,k))_{12}, 
\end{equation}
where $F(t,k)$ is some function whose first and second columns are analytic and bounded in $\C^-$ and $\C^+$ respectively and $F(t,k) = I + O(1/k)$ as $k \to \infty$ with $k \in (\C^-, \C^+)$.
We can write the $(12)$ entry of the global relation (\ref{globalrelation}) as 
\begin{align}\label{GRc}
& c(t, k) = \Phi_1(t, k) + \frac{b(k)}{a(k)}\overline{\Phi_2(t, \bar{k})} e^{-2if_2(k)t}, \qquad \text{Im}\,k \geq 0.
\end{align} 
The function $c(t, k)$ is analytic and bounded in $\text{Im}\,k > 0$ away from the possible zeros of $a(k)$. The functions $a(k)$ and $b(k)$ are analytic and bounded in $D_1 \cup D_2$, while $A(k)$ and $B(k)$ are entire functions which are bounded in $D_1 \cup D_3$. The functions $\Phi_1(t, k)$ and $\Phi_2(t, k)$ are entire functions which are bounded forÊ $k \in D_2 \cup D_4$. 
Furthermore, the second column vector of equation (\ref{Phivolterra}) yields
\begin{align}\nonumber
& \Phi_1(t,k) = \int_0^t e^{-4ik^2(t - \tau)}\left[(\tilde{Q}_0)_{11} \Phi_1 + (\tilde{Q}_0)_{12} \Phi_2\right](\tau, k) d\tau, 
  	\\ \label{Phi1Phi2system}
& \Phi_2(t,k) = 1 +  \int_0^t \left[(\tilde{Q}_0)_{21} \Phi_1 + (\tilde{Q}_0)_{22} \Phi_2\right](\tau, k) d\tau, 	 \qquad 0 < t < T, \; k \in \C.
\end{align}
Hence,
\begin{subequations}\label{hatPhieqs}
\begin{align}
& \hat{\Phi}_1(t,k) = \int_0^t \left[(\tilde{Q}_0)_{11} \hat{\Phi}_1 + (\tilde{Q}_0)_{12} \hat{\Phi}_2\right](\tau, k) d\tau, 
  	\\
& \hat{\Phi}_2(t,k) = e^{4ik^2t} +  \int_0^t e^{4ik^2(t - \tau)} \left[(\tilde{Q}_0)_{21} \hat{\Phi}_1 + (\tilde{Q}_0)_{22} \hat{\Phi}_2\right](\tau, k) d\tau, 	 \qquad 0 < t < T, \; k \in \C,
\end{align}
\end{subequations}
where
\begin{align}\label{hatPhiPhie}
\hat{\Phi}_j(t,k) = \Phi_j(t,k) e^{4ik^2t}, \qquad j = 1,2.
\end{align}
Therefore, $\hat{\Phi}_j(t,k)$, $j = 1,2$, is bounded for $k \in D_1 \cup D_3$ and $\overline{\Phi_j(t,\bar{k})}e^{-4ik^2t}$ is bounded for $k \in D_2 \cup D_4$.
Equation (\ref{GRc}) shows that whereas $\Phi_1$ and $\bar{\Phi}_2e^{-4ik^2t}$ are bounded in $D_2 \cup D_4$, the combination appearing on the right-hand side of (\ref{GRc}) is also bounded in $D_1$.

\subsection{Asymptotics}
Integration by parts in (\ref{Phi1Phi2system}) shows that
\begin{align}\nonumber
& \Phi_1(t, k) = \frac{\Phi_1^{(1)}(t)}{k} + \frac{\Phi_1^{(2)}(t)}{k^2} + O \Bigl(\frac{1}{k^3} \Bigr) + O\Bigl(\frac{e^{-4ik^2t}}{k}\Bigr), \qquad k \to \infty, \quad k \in D_2 \cup D_4,
	\\\label{cPhiexpansions}
& \Phi_2(t, k) = 1 + \frac{\Phi_2^{(1)}(t)}{k}  + O \Bigl(\frac{1}{k^2} \Bigr), \qquad k \to \infty, \quad k \in D_2 \cup D_4,
\end{align}
where
\begin{align*}
& \Phi_1^{(1)}(t) = \frac{g_0(t)}{2i}, \qquad
\Phi_1^{(2)}(t) = \frac{g_1(t)}{4} - \frac{ig_0(t)}{2}\int_{(0, 0)}^{(0,t)} \Delta(\xi, \tau), 
	\\
& \Phi_2^{(1)}(t) = \int_{(0, 0)}^{(0,t)} \Delta(\xi, \tau),
\end{align*}
and the closed one-form $\Delta$ is defined by
$$\Delta(x, t) = \frac{\lambda}{2}\left[-i|q(x, t)|^2dx + (\bar{q}(x, t)q_x(x, t) - q(x, t)\bar{q}_x(x, t))dt\right].$$
In particular, we find the following expressions for the boundary values:
\begin{subequations}
\begin{align}\label{g0Phi11}
& g_0(t) = 2i \Phi_1^{(1)}(t),
	\\ \label{g1cPhi}
& g_1(t) = 4 \Phi_1^{(2)}(t) + 2ig_0\Phi_2^{(1)}(t), \qquad 0 \leq t < T.
\end{align}
\end{subequations}

We will also need the asymptotics of $c$.
\begin{lemma}\label{clemma}
The global relation (\ref{GRc}) implies that the large $k$ behavior of $c(t,k)$ satisfies
\begin{align}\label{casymptotics} 
 c(t, k) = \frac{\Phi_1^{(1)}(t)}{k} + \frac{\Phi_1^{(2)}(t)}{k^2} + O \Bigl(\frac{1}{k^3} \Bigr),
 \qquad k \to \infty, \quad k \in D_1.
\end{align}
\end{lemma}
\proofbegin
See appendix \ref{clemmaapp}.
\proofend

\subsection{The Dirichlet and Neumann problems}
The following theorem expresses the spectral functions $A(k)$ and $B(k)$ in terms of the prescribed initial and boundary data via the solution of a system of nonlinear integral equations. This result was already derived in \cite{F2005} using the GLM representations for $\Phi_1$ and $\Phi_2$. However, the derivation here does not require the use of the GLM representations; the result is instead deduced directly from the asymptotics of the eigenfunctions as $k \to \infty$. 

We let $\partial D_j$, $j = 1, \dots, 4$, denote the boundary of the $j$'th quadrant $D_j$, oriented so that $D_j$ lies to the left of $\partial D_j$.

\begin{theorem}\label{th1}
Let $T < \infty$. Let $q_0(x)$, $x \geq 0$, be a function of Schwartz class. For the Dirichlet problem it is assumed that the function $g_0(t)$, $0 \leq t < T$, has sufficient smoothness and is compatible with $q_0(x)$ at $x=t=0$. Similarly, for the Neumann problem it is assumed that the function Ê$g_1(t)$, $0 \leq t < T$, has sufficient smoothness and is compatible with $q_0(x)$ at $x=t=0$.
Suppose that $a(k)$ has a finite (possibly empty) set of simple zeros, which are denoted by $\{k_j\}_1^N$; assume that no zeros occur on the boundaries of $D_1$ and $D_2$. 

Then the spectral functions $A(k)$ and $B(k)$ associated with the NLS equation (\ref{nls}) are given by
\begin{align}\label{ABexpressions}
& A(k)  = \overline{\Phi_2(T, \bar{k})} 
\qquad B(k) = -\Phi_1(T, k)e^{4ik^2T},
\end{align}
where the complex-valued functions $\Phi_1(t, k)$ and $\Phi_2(t, k)$ satisfy the following system of nonlinear integral equations:
\begin{align} \label{Phieqs}
&  \Phi_1(t, k) = \int_0^t e^{4ik^2(t' - t)} \bigl[-i\lambda |g_0|^2 \Phi_1 + (2kg_0 + ig_1)\Phi_2\bigr] (t', k) dt',
	\\ \nonumber
&   \Phi_2(t, k) = 1 + \lambda \int_0^t \bigl[(2k\bar{g}_0  - i\bar{g}_1)\Phi_1 + i|g_0|^2\Phi_2\bigr](t', k) dt', \qquad 0 < t < T, \; k \in \C.
\end{align}
\begin{itemize}
\item[$(a)$] For the Dirichlet problem, the unknown Neumann boundary value $g_1(t)$ is given by
\begin{align} \nonumber
g_1(t) = \;& \frac{2}{\pi i}\int_{\partial D_3}\bigl(k\chi_1(t,k) + i g_0(t)\bigr) dk + \frac{2g_0(t)}{\pi} \int_{\partial D_3} \chi_2(t,k) dk
	\\ \nonumber
& - \frac{4}{\pi i}\int_{\partial D_3} ke^{-4ik^2 t}\frac{b(-k)}{a(-k)} \overline{\Phi_2(t, -\bar{k})}  dk
	\\ \label{g1expression}
& + 8 \sum_{k_j \in D_1} k_je^{-4ik_j^2 t}\frac{b(k_j)}{\dot{a}(k_j)} \overline{\Phi_2(t, \bar{k}_j)}, \qquad 0 < t < T,
\end{align}
where $\{\chi_j\}_1^2$ denote the odd combinations formed from Ê$\{\Phi_j\}_1^2$, i.e. 
\begin{equation}\label{chijdef}
  \chi_j(t,k) = \Phi_j(t,k) - \Phi_j(t, -k), \qquad j = 1,2, \quad 0 < t < T, \quad k \in \C.
\end{equation}  

\item[$(b)$] For the Neumann problem, the unknown Dirichlet boundary value $g_0(t)$ is given by
\begin{align}\label{g0expression}  
  g_0(t) =\;& \frac{1}{\pi}\int_{\partial D_3} \check{\chi}_1(t,k) dk +  \frac{2}{\pi} \int_{\partial D_3} e^{-4ik^2t} \frac{b(-k)}{a(-k)}\overline{\Phi_2(t, -\bar{k})} dk
  	\\ \nonumber
&   + 4i\sum_{k_j \in D_1} e^{-4ik_j^2 t}\frac{b(k_j)}{\dot{a}(k_j)} \overline{\Phi_2(t, \bar{k}_j)},
\end{align}
where $\{\check{\chi}_j\}_1^2$ denote the even combinations formed from Ê$\{\Phi_j\}_1^2$, i.e.
\begin{align}\label{checkchijdef}
\check{\chi}_j(t,k) = \Phi_j(t,k) + \Phi_j(t,-k), \qquad  j = 1,2, \quad 0 < t < T, \quad k \in \C.
\end{align}
\end{itemize}
\end{theorem}
\proofbegin
Equations (\ref{ABexpressions}) and (\ref{Phieqs}) follow from the definition of $\Phi_1$ and $\Phi_2$ and from (\ref{Phivolterra}).


$(a)$ In order to derive (\ref{g1expression}) we note that equation (\ref{g1cPhi}) expresses $g_1$ in terms of $\Phi_2^{(1)}$ and $\Phi_1^{(2)}$. Furthermore, equations (\ref{cPhiexpansions}) and Cauchy's theorem imply
\begin{align}\label{Phi21fromPhi2}
  -\frac{i\pi}{2} \Phi_2^{(1)}(t) = \int_{\partial D_2} [\Phi_2(t,k) -1] dk = \int_{\partial D_4} [\Phi_2(t,k) - 1] dk
\end{align}
and
\begin{align}\label{Phi12fromPhi1}
  -\frac{i\pi}{2} \Phi_1^{(2)}(t) = \int_{\partial D_2} \biggl[k\Phi_1(t,k) - \frac{g_0(t)}{2i}\biggr] dk 
  = \int_{\partial D_4} \biggl[k\Phi_1(t,k) - \frac{g_0(t)}{2i}\biggr] dk,
\end{align}
where we have used equation (\ref{g0Phi11}) to write $\Phi_1^{(1)}$ in terms of $g_0(t)$.

If, instead of equation (\ref{g1expression}), we use equation (\ref{g1cPhi}) with $\Phi_2^{(1)}$ and $\Phi_1^{(1)}$ defined by (\ref{Phi21fromPhi2}) and (\ref{Phi12fromPhi1}), we do {\it not} obtain an effective characterization. The latter representation requires the appearance of $\chi_1$ and $\chi_2$. In this respect we note that 
\begin{align}\nonumber
  i\pi \Phi_2^{(1)}(t) & = -\biggl(\int_{\partial D_2} + \int_{\partial D_4}\biggr) [\Phi_2(t,k) -1] dk
  = \biggl(\int_{\partial D_3} + \int_{\partial D_1}\biggr) [\Phi_2(t,k) -1] dk
  	\\ \label{ipiPhi21}
&  = \int_{\partial D_3} [\Phi_2(t,k) -1] dk - \int_{\partial D_3} [\Phi_2(t,-k) -1] dk 
    = \int_{\partial D_3} \chi_2(t,k) dk.
\end{align}
Similarly,
\begin{align}\nonumber
  i\pi \Phi_1^{(2)}(t) & = \biggl(\int_{\partial D_3} + \int_{\partial D_1}\biggr) \biggl[k\Phi_1(t,k) - \frac{g_0(t)}{2i}\biggr] dk
  	\\ \nonumber
&  =  \biggl(\int_{\partial D_3} - \int_{\partial D_1}\biggr) \biggl[k\Phi_1(t,k) - \frac{g_0(t)}{2i}\biggr] dk + 2\int_{\partial D_1} \biggl[k\Phi_1(t,k) - \frac{g_0(t)}{2i}\biggr] dk
	\\ \label{ipiPhi12}
&  = \int_{\partial D_3}[k \chi_1(t,k) + ig_0(t)] dk + 2\int_{\partial D_1} \biggl[k\Phi_1(t,k) - \frac{g_0(t)}{2i}\biggr] dk.
\end{align}
The last step involves using the global relation (\ref{GRc}) in order to compute the second term on the right-hand side of (\ref{ipiPhi12}):
\begin{align}\nonumber
 2 \int_{\partial D_1} \biggl[k\Phi_1(t,k) - \frac{g_0(t)}{2i}\biggr] dk
  = &\; 2 \int_{\partial D_1} \biggl[kc(t,k) - \frac{g_0(t)}{2i}\biggr] dk
  	\\\label{laststep}
&  - 2\int_{\partial D_1} \frac{kb(k)}{a(k)}\overline{\Phi_2(t, \bar{k})} e^{-2if_2(k)t} dk.
\end{align}
Using the asymptotics (\ref{casymptotics}) of $c(t,k)$ and Cauchy's theorem to compute the first term on the right-hand side of (\ref{laststep}) and using the transformation $k \to -k$ in the second term on the right-hand side of (\ref{laststep}), we find
\begin{align} \nonumber
 2 \int_{\partial D_1} \biggl[k \Phi_1(t,k) - \frac{g_0(t)}{2i}\biggr] dk = &\; - i\pi \Phi_1^{(2)}(t) - 2 \int_{\partial D_3} \frac{kb(-k)}{a(-k)}\overline{\Phi_2(t, -\bar{k})}e^{-4ik^2t} 
  	\\ \label{secondtermcomputed}
& + 4\pi i \sum_{k_j \in D_1} \frac{k_jb(k_j)}{a(k_j)} \overline{\Phi_2(t, \bar{k}_j)}e^{-4ik_j^2t}.
\end{align}
Equations (\ref{ipiPhi21}), (\ref{ipiPhi12}), and (\ref{secondtermcomputed}) together with (\ref{g1cPhi}) yield (\ref{g1expression}).

$(b)$ In order to derive (\ref{g0expression}) we note that equation (\ref{g0Phi11}) expresses $g_0$ in terms of $\Phi_1^{(1)}$. Furthermore, equation (\ref{cPhiexpansions}) and Cauchy's theorem imply
\begin{align}
  -\frac{i\pi}{2} \Phi_1^{(1)}(t) = \int_{\partial D_2} \Phi_1(t,k) dk = \int_{\partial D_4} \Phi_1(t,k) dk.
\end{align}
Thus,
\begin{align}\nonumber
  i\pi \Phi_1^{(1)}(t) & = \biggl(\int_{\partial D_3} + \int_{\partial D_1}\biggr) \Phi_1(t,k) dk
  	\\ \nonumber
&  =  \biggl(\int_{\partial D_3} - \int_{\partial D_1}\biggr) \Phi_1(t,k)dk + 2\int_{\partial D_1} \Phi_1(t,k)dk
	\\ \label{ipiPhi11}
&  = \int_{\partial D_3} \check{\chi}_1(t,k) dk + 2\int_{\partial D_1} \Phi_1(t,k) dk.
\end{align}
The last step involves using the global relation (\ref{GRc}) and  to compute the second term on the right-hand side of (\ref{ipiPhi11}):
\begin{align}\label{Neumannlaststep}
2 \int_{\partial D_1} \Phi_1(t,k) dk
 & =  2\int_{\partial D_1} c(t,k) dk  - 2 \int_{\partial D_1} \frac{b(k)}{a(k)}\overline{\Phi_2(t, \bar{k})} e^{-2if_2(k)t} dk.
\end{align}
Using the asymptotics (\ref{casymptotics}) of $c(t,k)$ and Cauchy's theorem to compute the first term on the right-hand side of (\ref{Neumannlaststep}) and employing the transformation $k \to -k$ in the second term on the right-hand side of (\ref{Neumannlaststep}), we find
\begin{align} \nonumber
2 \int_{\partial D_1} \Phi_1(t,k) dk
 = &\; - i\pi \Phi_1^{(1)}(t) + 2\int_{\partial D_3} \frac{b(-k)}{a(-k)}\overline{\Phi_2(t, -\bar{k})}e^{-4ik^2t} 
  	\\ \label{Neumannsecondtermcomputed}
& + 4\pi i \sum_{k_j \in D_1} \frac{b(k_j)}{a(k_j)} \overline{\Phi_2(t, \bar{k}_j)}e^{-4ik_j^2t}.
\end{align}
Equations (\ref{ipiPhi11}) and (\ref{Neumannsecondtermcomputed}) together with (\ref{g0Phi11}) yield (\ref{g0expression}).
\proofend


\begin{remark}\upshape
The Dirichlet and Neumann problems for the NLS equation on the half-line can now be solved as follows: If $a(k)$ has no zeros, the functions $A(k)$ and $B(k)$ enter the formulation of the RH problem of subsection \ref{RHsubsec} only in the combination $A(k)/B(k)$ for $k \in \partial D_3$. 
For the Dirichlet problem, substituting the expression (\ref{g1expression}) for $g_1(t)$ into (\ref{Phieqs}) yields a system of quadratically nonlinear integral equations involving the functions $\Phi_1(t, k)$ and $\Phi_2(t,k)$, $0 < t < T$, $k \in \R \cup i\R$. For the Neumann problem, substituting the expression (\ref{g0expression}) for $g_0(t)$ into (\ref{Phieqs}) yields a system of quadratically nonlinear integral equations for $\Phi_1(t, k)$ and $\Phi_2(t,k)$.
Assuming that these systems have unique solutions, $A(k)/B(k)$, $k \in \partial D_3$, can be determined from (\ref{ABexpressions}). 
In fact, we will show in the following subsection that these nonlinear systems provide {\it effective} characterizations of $A(k)$ and $B(k)$ in the sense defined in the introduction. In particular, the systems can be solved recursively to all orders in a perturbative scheme. If $a(k)$ has zeros, the sum of residues also needs to be taken into account.
\end{remark}

\subsection{Effective characterization}\label{NLSeffectivesubsec}
Substituting into the system (\ref{Phieqs}) the expansions
\begin{align}
 \Phi_j &= \Phi_{j0} + \epsilon \Phi_{j1} + \epsilon^2 \Phi_{j2} + \cdots, \qquad \epsilon \to 0, \; j = 1,2,
  	\\\label{g0g1expansions}
 g_0 &= \epsilon g_{01} + \epsilon^2 g_{02} + \cdots, \qquad  g_1 = \epsilon g_{11} + \epsilon^2 g_{12} + \cdots, 
 	\\
 h_0 &= \epsilon h_{01} + \epsilon^2 h_{02} + \cdots, \qquad  h_1 = \epsilon h_{11} + \epsilon^2 h_{12} + \cdots, 
\end{align}
where $\epsilon > 0$ is a small parameter, we find that the terms of $O(1)$ give $\Phi_{10} \equiv 0$ and $\Phi_{20} \equiv 1$; the terms of $O(\epsilon)$ give $\Phi_{21} \equiv 0$ and
\begin{align}\label{eps1a}
&O(\epsilon): \;  \Phi_{11}(t, k) = \int_0^t e^{4ik^2(t' - t)}  (2kg_{01}(t') + ig_{11}(t')) dt';
\end{align}
the terms of $O(\epsilon^2)$ give
\begin{subequations}\label{eps2}
\begin{align} \label{eps2a}
&O(\epsilon^2): \;   \Phi_{12}(t, k) = \int_0^t e^{4ik^2(t' - t)} (2kg_{02}(t') + ig_{12}(t')) dt';
	\\ \nonumber
&O(\epsilon^2): \;     \Phi_{22}(t, k) = \lambda \int_0^t \bigl[ (2k\bar{g}_{01}  - i\bar{g}_{11})\Phi_{11} +  i|g_{01}|^2\bigr](t', k) dt'.
\end{align}
\end{subequations}

On the other hand, expanding (\ref{g1expression}) and assuming for simplicity that $a(k)$ has no zeros, we find
\begin{align} \label{g11expression}
g_{11} = \;& \frac{2}{\pi i}\int_{\partial D_3}\bigl(k\chi_{11}(t,k) + i g_{01}(t)\bigr) dk 
- \frac{4}{\pi i}\int_{\partial D_3} ke^{-4ik^2 t} b_1(-k)dk,
\end{align}
where $\chi_1 = \epsilon \chi_{11} + O(\epsilon^2)$ and $b = \epsilon b_1 + O(\epsilon^2)$. 

The Dirichlet problem can now be solved perturbatively as follows: Let $\chi_j$, $j=1,2$, denote the odd combination in $k$ formed from $\Phi_j$ as defined in (\ref{chijdef}).
The odd part of (\ref{eps1a}) yields
\begin{align}\label{chi11eqn}
& \chi_{11}(t, k) = 4k \int_0^t e^{4ik^2(t' - t)} g_{01}(t') dt'.
\end{align}
Given $g_{01}$, we can use this equation to determine $\chi_{11}(t,k)$. We can then compute $g_{11}$ from (\ref{g11expression}) and then $\Phi_{11}$ follows from (\ref{eps1a}). 
In the same way, we can use the odd part of (\ref{eps2a}) to determine $\chi_{12}$; we use $\chi_{12}$ to compute $g_{12}$, and then $\Phi_{12}$ and $\Phi_{22}$ follow from (\ref{eps2}). This recursive scheme can be continued indefinitely. Indeed, suppose $\Phi_{1j}$, $\Phi_{1j}$, and $g_{1j}$ have been determined for all $0 \leq j \leq n-1$, for some $n \geq 0$. The terms in (\ref{Phieqs}) of $O(\epsilon^n)$ give
\begin{subequations}\label{epsn}
\begin{align} \label{epsna}
&O(\epsilon^n): \;  
\Phi_{1n}(t, k) = \int_0^t e^{4ik^2(t' - t)}  (2kg_{0n}(t') + ig_{1n}(t')) dt' + \lot,
	\\ \label{epsnb}
 &O(\epsilon^n): \;   \Phi_{2n} = \lot,
\end{align}
\end{subequations}
where `$\lot$' denotes an expression involving known terms of lower order. Similarly, the terms of $O(\epsilon^n)$ of the integral representation (\ref{g1expression}) for $g_1$ give 
\begin{equation}\label{g1nexpression}
g_{1n} = \frac{2}{\pi i}\int_{\partial D_3}\bigl(k\chi_{1n} + i g_{0n}\bigr) dk
- \frac{4}{\pi i}\int_{\partial D_3} ke^{-4ik^2 t}b_n(-k) dk + \lot.
\end{equation}
The odd part of (\ref{epsna}) yields 
$$\chi_{1n}(t, k) = 4k \int_0^t e^{4ik^2(t' - t)} g_{0n}(t') dt' + \lot.$$
Substituting the solution $\chi_{1n}$ of this equation into (\ref{g1nexpression}), we find $g_{1n}$ and then $\Phi_{1n}$ and $\Phi_{2n}$ are found from (\ref{epsn}).
This shows that for the Dirichlet problem $\Phi_1$ and $\Phi_2$ can be determined to all orders in a perturbative scheme by solving the nonlinear system of theorem \ref{th1} recursively.

Similarly, for the Neumann problem, substituting the expression (\ref{g0expression}) for $g_0(t)$ into (\ref{Phieqs}) yields a system of quadratically nonlinear integral equations for the functions $\Phi_1(t, k)$ and $\Phi_2(t,k)$. This nonlinear system can be solved recursively to all orders perturbatively and therefore it provides an effective characterization of $A(k)$ and $B(k)$ for the Neumann problem. Indeed, letting
$$\check{\chi}_1 = \epsilon \check{\chi}_{11} +  \epsilon \check{\chi}_{12} + \cdots,$$
the even part of (\ref{eps1a}) yields
\begin{align}\label{checkchi11t}
  \check{\chi}_{11}(t, k) = 2i \int_0^t e^{4ik^2(t' - t)} g_{11}(t') dt',
\end{align}
while (\ref{g0expression}) yields
\begin{equation}\label{g01expression}
 g_{01} = \frac{1}{\pi}\int_{\partial D_3} \check{\chi}_{11}(t,k) dk +  \frac{2}{\pi} \int_{\partial D_3} e^{-4ik^2t} b_1(-k) dk.
\end{equation}
Since $g_1$ is known, (\ref{checkchi11t}) can be solved for $\check{\chi}_{11}$, and then $g_{01}$ can be determined from (\ref{g01expression}). After $g_{01}$ has been found, $\Phi_{11}$ can be computed from (\ref{eps1a}). Extending this procedure to higher orders we find, just like in the case of the Dirichlet problem, a recursive scheme which can be used to determine $\Phi_1$ and $\Phi_2$ to all orders.

\begin{remark}\upshape
The function $g_0$, in addition to (\ref{g0expression}), also admits the following alternative representation:
\begin{equation}\label{alternativeg0expression}  
  g_0(t) = \frac{2}{\pi}\int_{\partial D_3} \chi_1(t,k) dk.
\end{equation}
However, this alternative representation is not suitable for the effective solution of the Neumann problem. Indeed, for the Neumann problem, we can find $\check{\chi}_{1n}$ from the terms of $O(\epsilon^n)$, so that the representation (\ref{g0expression}) can be used to find $g_{0n}$. However, the function $\chi_{1n}$ remains unknown, so that the representation (\ref{alternativeg0expression}) cannot be used to find $g_{0n}$. 
\end{remark}

\section{The mKdV equation}\label{mKdVIsec}\nequation
We consider the mKdV equation (\ref{mkdvI}) posed on the half-line.

\subsection{Asymptotics}

In this case, the $\Phi_j$'s admit the asymptotics
\begin{subequations} \label{Phiabcexpansions}
\begin{align} \label{Phiabcexpansionsa}
& \Phi_1(t, k) = \frac{\Phi_1^{(1)}(t)}{k} + \frac{\Phi_1^{(2)}(t)}{k^2} + \frac{\Phi_1^{(3)}(t)}{k^3}+ O \Bigl(\frac{1}{k^4} \Bigr) + O\Bigl(\frac{e^{-8ik^3t}}{k^2}\Bigr), 
	\\ \label{Phiabcexpansionsb}
& \Phi_2(t, k) = 1 + \frac{\Phi_2^{(1)}(t)}{k} + \frac{\Phi_2^{(2)}(t)}{k^2}  + O \Bigl(\frac{1}{k^3} \Bigr), \qquad k \to \infty, \quad k \in D_2 \cup D_4,
\end{align}
\end{subequations}
where 
\begin{align*}
& 
\Phi_2^{(1)} = \frac{\lambda}{2i} \int_{(0, 0)}^{(0,t)} \Delta, \qquad
\Phi_2^{(2)} = \int_{(0, 0)}^{(0,t)} \tilde{\Delta}, \qquad 
\Phi_1^{(1)} = \frac{g_0}{2i}, 
	\\
& \Phi_1^{(2)} = \frac{g_1}{4} - \frac{\lambda g_0}{4}\int_{(0, 0)}^{(0,t)} \Delta, \qquad \Phi_1^{(3)} = \frac{1}{2i}\biggl(g_0\Phi_2^{(2)} + \frac{i}{2}g_1\Phi_2^{(1)} + \frac{\lambda g_0^3}{4} - \frac{1}{4}g_2\biggr),
\end{align*}
and the closed one-forms $\Delta$ and $\tilde{\Delta}$ are defined by
\begin{align*}
 \Delta & = q^2dx + (q_x^2 - 2qq_{xx} + 3\lambda q^4)dt,
\qquad
 \tilde{\Delta} 
 = \frac{1}{8}d\left[\lambda q^2 -  \left(\int_{(0,0)}^{(x,t)} \Delta\right)^2\right].
\end{align*}
In particular, we find the following expressions for the boundary values:
\begin{subequations}\label{mkdvg0g1g2}
\begin{align}\label{mkdvg0Phi}
& g_0 = 2i\Phi_1^{(1)},
    	\\\label{mkdvg1Phi}
& g_1 = 2ig_0\Phi_2^{(1)} + 4 \Phi_1^{(2)},
	\\\label{mkdvg2Phi}
& g_2 = \lambda g_0^3 - 8i\Phi_{1}^{(3)} + 4g_0\Phi_{2}^{(2)} + 2ig_1\Phi_{2}^{(1)}.
\end{align}
\end{subequations}

We will also need the asymptotics of the function $c(t,k)$ defined in (\ref{cddef}).
\begin{lemma}\label{mKdVclemma}
The global relation (\ref{GRc}) implies that the large $k$ behavior of $c(t,k)$ satisfies
\begin{align} \label{mkdvcasymptotics}
& c(t, k) = \frac{\Phi_1^{(1)}(t)}{k} + \frac{\Phi_1^{(2)}(t)}{k^2} + \frac{\Phi_1^{(3)}(t)}{k^3} +  O \Bigl(\frac{1}{k^4} \Bigr), 
\qquad k \to \infty, \; k \in D_1.
\end{align}
\end{lemma}
\proofbegin See appendix \ref{clemmaapp}. \proofend

\subsection{The Dirichlet and Neumann problems}
We can now derive effective characterizations of $A(k)$ and $B(k)$ for the Dirichlet ($g_0$ prescribed), the first Neumann ($g_1$ prescribed), and the second Neumann ($g_2$ prescribed) problems. 

Define $\alpha$ by $\alpha = e^{2\pi i/3}$ and let $\{\chi_j, \hat{\chi}_j, \check{\chi}_j\}_1^2$ denote the following combinations formed from $\{\Phi_j\}_1^2$: 
\begin{align*}
& \chi_j(t,k) = \Phi_j(t, k) + \alpha\Phi_j(t, \alpha k)+ \alpha^2 \Phi_j(t, \alpha^2 k), \qquad j = 1,2,
	\\	
& \hat{\chi}_j(t,k) = \Phi_j(t, k) + \alpha^2 \Phi_j(t, \alpha k)+ \alpha \Phi_j(t, \alpha^2 k), \qquad j = 1,2,
	\\
& \check{\chi}_j(t,k) = \Phi_j(t,k) + \Phi_j(t,\alpha k) + \Phi_j(t, \alpha^2 k), \qquad j = 1,2.
\end{align*}
Define $R(t, k)$ by\footnote{The $t$-dependence of $R(t,k)$ will be suppressed below.}
$$R(t, k) = \frac{b(k)\overline{\Phi_2(t, \bar{k})}}{a(k)}, \qquad 0 < t < T, \; k \in \C.$$
Let $D_1 = D_1' \cup D_1''$ where $D_1' = D_1 \cap \{\text{Re}\, k > 0\}$ and $D_1'' = D_1 \cap \{\text{Re}\, k < 0\}$. Similarly, let $D_4 = D_4' \cup D_4''$ with $D_4' = D_4 \cap \{\text{Re}\, k > 0\}$ and $D_4'' = D_4 \cap \{\text{Re}\, k < 0\}$.

\begin{theorem}\label{mkdvth1}
Let $T < \infty$. Let $q_0(x)$, $x \geq 0$, be a function of Schwartz class. For the Dirichlet problem it is assumed that the function $g_0(t)$, $0 \leq t < T$, has sufficient smoothness and is compatible with $q_0(x)$ at $x=t=0$. Similarly, for the first and second Neumann problems it is assumed that the functions Ê$g_1(t)$ and $g_2(t)$, $0 \leq t < T$, have sufficient smoothness and are compatible with $q_0(x)$ at $x=t=0$, respectively.
Suppose that $a(k)$ has a finite (possibly empty) set of simple zeros, which are denoted by $\{k_j\}_1^N$; assume that no zeros occur on the boundaries of $D_1$ and $D_2$. 

Then the spectral functions $A(k)$ and $B(k)$ associated with the mKdV equation (\ref{mkdvI}) are given by
\begin{align}\label{mkdvABexpressions}
& A(k)  = \overline{\Phi_2(T, \bar{k})} 
\qquad B(k) = -\Phi_1(T, k)e^{8ik^3T},
\end{align}
where the complex-valued functions $\Phi_1(t, k)$ and $\Phi_2(t, k)$ satisfy the following system of nonlinear integral equations:
\begin{subequations}\label{mkdvPhieqs}
\begin{align} \label{mkdvPhieqsa}
&  \Phi_1(t, k) = \int_0^t e^{8ik^3(t' - t)} \bigl[-2 i k \lambda  g_0^2 \Phi_1 + \left( 2 \lambda  g_0^3+4 k^2 g_0+2 i k g_1 - g_2\right)\Phi_2\bigr] (t', k) dt',
	\\ \nonumber
&   \Phi_2(t, k) = 1 + \lambda \int_0^t \bigl[ \left(2 \lambda  g_0^3+4 k^2 g_0-2 i k g_1 - g_2\right)\Phi_1 + 2 i k  g_0^2\Phi_2\bigr](t', k) dt', 
	\\ \label{mkdvPhieqsb}
&\hspace{8cm} 0 < t < T, \; k \in \C.
\end{align}
\end{subequations}

\begin{itemize}
\item[$(a)$] For the Dirichlet problem, the unknown Neumann boundary values $g_1(t)$ and $g_2(t)$ are given by
\begin{subequations}\label{mkdvg1g2expression}
\begin{align} \nonumber
g_1(t) = \;& \frac{2g_0(t)}{\pi}\int_{\partial D_3} \chi_2(t,k) dk 
+ \frac{2}{\pi i}\int_{\partial D_3} \left[k\chi_1(t,k) - \frac{3g_0(t)}{2i}\right] dk
	\\ \nonumber
& 
- \frac{2}{\pi i} \int_{\partial D_3} k e^{-8ik^3t} \bigl[(\alpha^2 - \alpha) R(\alpha k) + (\alpha - \alpha^2)R(\alpha^2 k) \bigr] dk
  	\\ \label{mkdvg1expression}
& + 4\biggl\{(1-\alpha^2)\sum_{k_j \in D_1'} + (1-\alpha)\sum_{k_j \in D_1''}\biggr\} k_j e^{-8ik_j^3t} \underset{k_j}{\res}R(k)
\end{align}
and
\begin{align} \nonumber
g_2(t) = \;& \lambda g_0^3(t) - \frac{4}{\pi} \int_{\partial D_3} \biggl[k^2\chi_1(t,k) - \frac{3kg_0(t)}{2i}\biggr] dk
 	\\ \nonumber
&  + \frac{4}{\pi} \int_{\partial D_3} 
k^2 e^{-8ik^3 t} \bigl[(1-\alpha)R(\alpha k) + (1 - \alpha^2) R(\alpha^2 k)\bigr]dk
	\\ \nonumber
& - 8i \biggl\{(1 - \alpha) \sum_{k_j \in D_1'} + (1 - \alpha^2)\sum_{k_j \in D_1''}\biggr\} k_j^2 e^{-8ik_j^3t} \underset{k_j}{\res}R(k)
	\\ \label{mkdvg2expression}
& + \frac{4g_0(t)}{\pi i}\int_{\partial D_3} k \hat{\chi}_2(t,k) dk + \frac{2g_1(t)}{\pi}\int_{\partial D_3} \chi_2(t,k) dk.
\end{align}
\end{subequations}

\item[$(b)$] For the first Neumann problem, the unknown boundary values $g_0(t)$ and $g_2(t)$ are given by
\begin{subequations}\label{N1expressions}
\begin{align} \nonumber
  g_0(t) =\;& \frac{1}{\pi} \int_{\partial D_3} \hat{\chi}_1(t,k) dk 
- \frac{1}{\pi} \int_{\partial D_3} e^{-8ik^3t} \bigl[(\alpha - \alpha^2)R(\alpha k) + (\alpha^2 - \alpha)R(\alpha^2 k)\bigr]dk
	\\ \label{N1g0expression} 
& + 2i\biggl\{(1 - \alpha) \sum_{k_j \in D_1'} + (1 - \alpha^2)\sum_{k_j \in D_1''}\biggr\}e^{-8ik_j^3t}\underset{k_j}{\res}R(k)
\end{align}
and
\begin{align} \nonumber
  g_2(t) =\;& \lambda g_0^3(t) -
 \frac{4}{\pi} \int_{\partial D_3} \left(k^2\hat{\chi}_1(t,k) - \frac{3}{\pi i}\int_{\partial D_3} l \hat{\chi}_1(t,l) dl\right) dk 
	\\ \nonumber
&  + \frac{4}{\pi} \int_{\partial D_3} 
k^2 e^{-8ik^3 t} \bigl[(1 - \alpha^2)R(\alpha k) + (1 - \alpha)R(\alpha^2 k)\bigr]dk
	\\ \nonumber
&- 8i\biggl\{(1- \alpha^2) \sum_{k_j \in D_1'} + (1 - \alpha)\sum_{k_j \in D_1''} \biggr\} k_j^2 e^{-8ik_j^3t}\underset{k_j}{\res} R(k)
	\\  \label{N1g2expression}  
& + \frac{4g_0(t)}{\pi i}\int_{\partial D_3} k \hat{\chi}_2(t,k) dk + \frac{2g_1(t)}{\pi}\int_{\partial D_3} \chi_2(t,k) dk.
\end{align}
\end{subequations}

\item[$(c)$] For the second Neumann problem, the unknown boundary values $g_0(t)$ and $g_1(t)$ are given by
\begin{subequations}\label{N2expressions}
\begin{align}\nonumber
  g_0(t) =\;&\frac{1}{\pi} \int_{\partial D_3} \check{\chi}_1(t,k) dk - \frac{1}{\pi} \int_{\partial D_3} e^{-8ik^3t} \bigl[(\alpha - 1)R(\alpha k) + (\alpha^2 - 1)R(\alpha^2 k)\bigr]dk
	\\ \label{N2g0expression}
& + 2i\biggl\{(1 - \alpha^2) \sum_{k_j \in D_1'} + (1 - \alpha)\sum_{k_j \in D_1''}\biggr\}e^{-8ik_j^3t}\underset{k_j}{\res}R(k)
\end{align}
and
\begin{align}\nonumber
  g_1(t) =\;& \frac{2g_0(t)}{\pi}\int_{\partial D_3} \chi_2(t,k)dk 
+   \frac{2}{\pi i} \int_{\partial D_3} k\check{\chi}_1(t,k) dk 
	\\ \nonumber
& - \frac{2}{\pi i} \int_{\partial D_3}k e^{-8ik^3t} \bigl[(\alpha^2 - 1) R(\alpha k) + (\alpha - 1)R(\alpha^2 k) \bigr]  dk
  	\\ \label{N2g1expression}
& + 4\biggl\{(1-\alpha)\sum_{k_j \in D_1'} + (1-\alpha^2)\sum_{k_j \in D_1''}\biggr\} k_j e^{-8ik_j^3t} \underset{k_j}{\res}R(k).
\end{align}
\end{subequations}
\end{itemize}
\end{theorem}
\proofbegin
$(a)$ In order to derive (\ref{mkdvg1expression}) we note that equation (\ref{mkdvg1Phi}) expresses $g_1$ in terms of $\Phi_2^{(1)}$ and $\Phi_1^{(2)}$. Furthermore, equations (\ref{Phiabcexpansions}) and Cauchy's theorem imply
\begin{align*}
  -\frac{2i\pi}{3} \Phi_2^{(1)}(t) = 2\int_{\partial D_2} [\Phi_2(t,k) -1] dk = \int_{\partial D_4} [\Phi_2(t,k) - 1] dk
\end{align*}
and
\begin{align*}
  -\frac{2i\pi}{3} \Phi_1^{(2)}(t) = 2\int_{\partial D_2} \biggl[k\Phi_1(t,k) - \frac{g_0(t)}{2i}\biggr] dk 
  = \int_{\partial D_4} \biggl[k\Phi_1(t,k) - \frac{g_0(t)}{2i}\biggr] dk.
\end{align*}
Thus,
\begin{align}\nonumber
  i\pi &\Phi_2^{(1)}(t) = -\biggl(\int_{\partial D_2} + \int_{\partial D_4}\biggr) [\Phi_2(t,k) -1] dk
  = \biggl(\int_{\partial D_3} + \int_{\partial D_1}\biggr) [\Phi_2(t,k) -1] dk
  	\\ \nonumber
&  = \int_{\partial D_3} [\Phi_2(t,k) -1] dk + \alpha \int_{\partial D_3} [\Phi_2(t,\alpha k) -1] dk 
+ \alpha^2 \int_{\partial D_3} [\Phi_2(t,\alpha^2 k) -1] dk 
	\\ \label{mkdvipiPhi21}
&   = \int_{\partial D_3} \chi_2(t,k) dk.
\end{align}
Similarly, 
\begin{align}\nonumber
  i\pi \Phi_1^{(2)}(t) = &\; \biggl(\int_{\partial D_3} + \int_{\partial D_1}\biggr) \biggl[k\Phi_1(t,k) - \frac{g_0(t)}{2i}\biggr] dk
  	\\ \nonumber
 =&\;  \biggl(\int_{\partial D_3} + \alpha^2 \int_{\partial D_1'} + \alpha \int_{\partial D_1''}\biggr) \biggl[k\Phi_1(t,k) - \frac{g_0(t)}{2i}\biggr] dk + I(t)
 	\\ \label{mkdvipiPhi12}
 = &\; \int_{\partial D_3} \biggl[k \chi_1(t,k) - \frac{3g_0(t)}{2i}\biggr] dk + I(t),
\end{align}
where $I(t)$ is defined by
$$I(t) =  \biggl((1- \alpha^2)\int_{\partial D_1'} + (1- \alpha)\int_{\partial D_1''}\biggr) \biggl[k\Phi_1(t,k) - \frac{g_0(t)}{2i}\biggr] dk.$$
The last step involves using the global relation (\ref{GRc}) to compute $I(t)$:
\begin{align}\nonumber
I(t)  = &\; 
\biggl((1- \alpha^2)\int_{\partial D_1'} + (1- \alpha)\int_{\partial D_1''}\biggr) \biggl[kc(t,k) - \frac{g_0(t)}{2i}\biggr] dk
	\\ \label{mkdvlaststep}
& - \biggl((1- \alpha^2)\int_{\partial D_1'} + (1- \alpha)\int_{\partial D_1''}\biggr) k e^{-2if_2(k)t} R(k)  dk.
\end{align}
Using the asymptotics (\ref{mkdvcasymptotics}) of $c(t,k)$ and Cauchy's theorem to compute the first term on the right-hand side of (\ref{mkdvlaststep}) and using the transformations $k \to \alpha k$ and $k \to \alpha^2 k$ in the second term on the right-hand side of (\ref{mkdvlaststep}), we find
\begin{align} \nonumber
 I(t) = &\; -i\pi \Phi_1^{(2)}(t) - \int_{\partial D_3} k e^{-8ik^3t}  \bigl[(\alpha^2 - \alpha) R(\alpha k) + (\alpha - \alpha^2)R(\alpha^2 k) \bigr] dk
  	\\ \label{mkdvsecondtermcomputed}
& + 2\pi i \biggl\{(1-\alpha^2)\sum_{k_j \in D_1'} + (1-\alpha)\sum_{k_j \in D_1''}\biggr\} k_j e^{-8ik_j^3t} \underset{k_j}{\res}R(k).
\end{align}
Equations (\ref{mkdvipiPhi12}) and (\ref{mkdvsecondtermcomputed}) imply
\begin{align*}
 \Phi_1^{(2)}(t) = &\; \frac{1}{ 2\pi i}\int_{\partial D_3} \biggl[k \chi_1(t,k) - \frac{3g_0(t)}{2i}\biggr] dk
 	\\
&  - \frac{1}{ 2\pi i} \int_{\partial D_3} k e^{-8ik^3t} \bigl[(\alpha^2 - \alpha) R(\alpha k) + (\alpha - \alpha^2)R(\alpha^2 k) \bigr]  dk
  	\\ 
& + \biggl\{(1-\alpha^2)\sum_{k_j \in D_1'} + (1-\alpha)\sum_{k_j \in D_1''}\biggr\} k_j e^{-8ik_j^3t} \underset{k_j}{\res}R(k).
\end{align*}
This equation together with (\ref{mkdvg1Phi}) and (\ref{mkdvipiPhi21}) yields (\ref{mkdvg1expression}).

In order to derive (\ref{mkdvg2expression}), we note that (\ref{mkdvg2Phi}) expresses $g_2$ in terms of $\Phi_1^{(3)}$, $\Phi_2^{(2)}$, and $\Phi_2^{(1)}$. Equation (\ref{mkdvg2expression}) follows from the expression (\ref{mkdvipiPhi21}) for $\Phi_2^{(1)}$ and the following formulas:
\begin{subequations}\label{Dproofformulas}
\begin{align}\label{Dproofformulasa}
  \Phi_2^{(2)} = \;& \frac{1}{\pi i}\int_{\partial D_3} k\hat{\chi}_2 dk,
	\\ \nonumber
  \Phi_1^{(3)} =\;& \frac{1}{2\pi i}\int_{\partial D_3} \left(k^2\chi_1 - \frac{3kg_0}{2i}\right) dk
 	\\ \nonumber
&  - \frac{1}{2\pi i} \int_{\partial D_3} 
k^2 e^{-8ik^3 t} \bigl[(1 - \alpha)R(\alpha k) + (1 - \alpha^2)R(\alpha^2 k)\bigr]dk
	\\ 
& + \biggl\{(1- \alpha) \sum_{k_j \in D_1'} + (1 - \alpha^2)\sum_{k_j \in D_1''} \biggr\} k_j^2 e^{-8ik_j^3t}\underset{k_j}{\res} R(k) .
\end{align}
\end{subequations}

$(b)$ In order to derive the representations (\ref{N1expressions}) relevant for the first Neumann problem, we use (\ref{mkdvg0g1g2}) together with (\ref{mkdvipiPhi21}), (\ref{Dproofformulasa}), and the following formulas:
\begin{subequations}\label{N1proofformulas}
\begin{align} \nonumber
  \Phi_1^{(1)} = \;& \frac{1}{2\pi i} \int_{\partial D_3} \hat{\chi}_1 dk 
- \frac{1}{2\pi i} \int_{\partial D_3} e^{-8ik^3t} \bigl[(\alpha - \alpha^2)R(\alpha k) + (\alpha^2 - \alpha)R(\alpha^2 k)\bigr]dk
	\\ \label{N1Phi11formula}
& + \biggl\{(1 - \alpha) \sum_{k_j \in D_1'} + (1 - \alpha^2)\sum_{k_j \in D_1''}\biggr\}e^{-8ik_j^3t}\underset{k_j}{\res}R(k),
  	\\
\Phi_1^{(2)} = \;&\frac{1}{\pi i}\int_{\partial D_3} k \hat{\chi}_1 dk,
	\\ \nonumber
  \Phi_1^{(3)} = \;& \frac{1}{2\pi i} \int_{\partial D_3} \left(k^2\hat{\chi}_1 - 3\Phi_1^{(2)}\right) dk 
	\\ \nonumber
&  - \frac{1}{2\pi i} \int_{\partial D_3} 
k^2 e^{-8ik^3 t} \bigl[(1 - \alpha^2)R(\alpha k) + (1 - \alpha)R(\alpha^2 k)\bigr]dk
	\\ \label{N1Phi13formula}
& + \biggl\{(1- \alpha^2) \sum_{k_j \in D_1'} + (1 - \alpha)\sum_{k_j \in D_1''} \biggr\} k_j^2 e^{-8ik_j^3t}\underset{k_j}{\res} R(k).
\end{align}
\end{subequations}

$(c)$ In order to derive the representations (\ref{N2expressions}) relevant for the second Neumann problem, we use (\ref{mkdvg0g1g2}) together with (\ref{mkdvipiPhi21}) and the following formulas:
\begin{subequations}\label{N2proofformulas}
\begin{align} \nonumber
  \Phi_1^{(1)} = \;&\frac{1}{2\pi i} \int_{\partial D_3} \check{\chi}_1 dk - \frac{1}{2\pi i} \int_{\partial D_3} e^{-8ik^3t} \bigl[(\alpha - 1)R(\alpha k) + (\alpha^2 - 1)R(\alpha^2 k)\bigr]dk
	\\ \label{N2Phi11formula}
& + \biggl\{(1 - \alpha^2) \sum_{k_j \in D_1'} + (1 - \alpha)\sum_{k_j \in D_1''}\biggr\}e^{-8ik_j^3t}\underset{k_j}{\res}R(k),
  	\\  \nonumber
  \Phi_1^{(2)} =\;& \frac{1}{2\pi i} \int_{\partial D_3} k\check{\chi}_1 dk 
 - \frac{1}{ 2\pi i} \int_{\partial D_3}k e^{-8ik^3t} \bigl[(\alpha^2 - 1) R(\alpha k) + (\alpha - 1)R(\alpha^2 k) \bigr]  dk
  	\\ \label{N2Phi12formula}
& + \biggl\{(1-\alpha)\sum_{k_j \in D_1'} + (1-\alpha^2)\sum_{k_j \in D_1''}\biggr\} k_j e^{-8ik_j^3t} \underset{k_j}{\res}R(k).
 \end{align}
\end{subequations}

The proofs of the formulas (\ref{Dproofformulas})-(\ref{N2proofformulas}) rely on arguments similar to those used in the proof of (\ref{mkdvg1expression}).
\proofend

\begin{figure}
\begin{center}
\begin{overpic}[width=.4\textwidth]{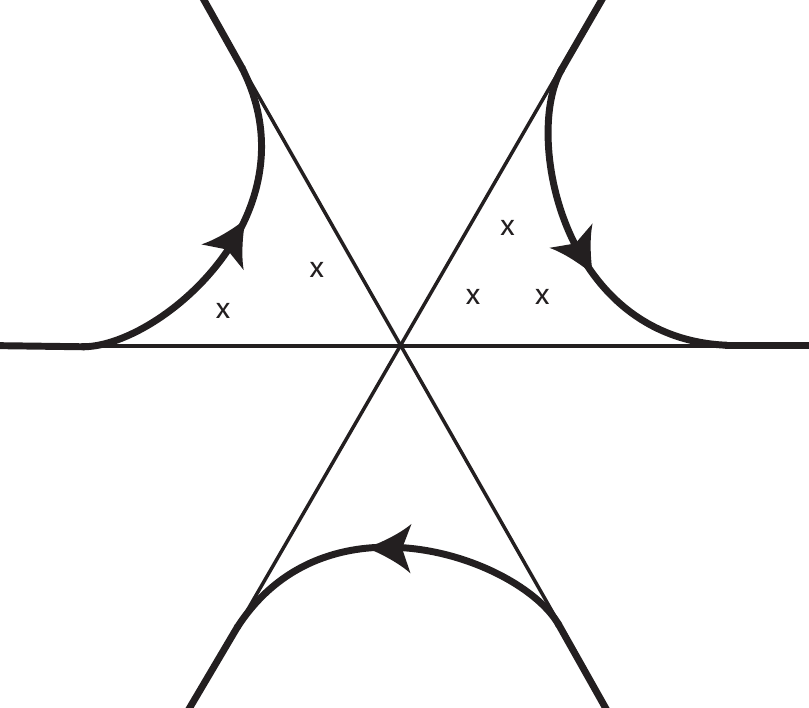}
      \put(90,70){$D_1'$}
      \put(0,70){$D_1''$}
      \put(47,10){$\Gamma$}
      \put(75,57){$\alpha \Gamma$}
      \put(13,57){$\alpha^2 \Gamma$}
       \put(103,43){$\text{Re}\; k$}
      \end{overpic}
     \begin{figuretext}\label{Gamma.pdf}
       The contour $\Gamma$ is defined in such a way that $\alpha \Gamma$ and $\alpha^2 \Gamma$ pass above the zeros of $a(k)$ (these zeros are indicated by $x$'s in the figure).
     \end{figuretext}
     \end{center}
\end{figure}

\begin{remark}\upshape
  The sums over the zeros of $a(k)$ that appear on the right-hand sides of (\ref{mkdvg1g2expression})-(\ref{N2expressions}) can be absorbed into the integrals involving $R(\alpha k)$ and $R(\alpha^2 k)$ by replacing the integration contour $\partial D_3$ for these integrals by $\Gamma$, where $\Gamma$ is a contour obtained by deforming $\partial D_3$ in such a way that $\alpha  \Gamma$ and $\alpha^2 \Gamma$ pass above all the zeros of $a(k)$ in $D_1'$ and $D_1''$, respectively, see figure \ref{Gamma.pdf}.
\end{remark}

\subsection{Effective characterization}
The nonlinear system for $\Phi_1$ andÊ $\Phi_2$ obtained by substituting the expressions (\ref{mkdvg1g2expression}) for $g_1$ and $g_2$ into (\ref{mkdvPhieqs}) provides an effective characterization for $A(k)$ and $B(k)$ for the Dirichlet problem for mKdV. 
Similarly, substituting the representations (\ref{N1expressions}) and (\ref{N2expressions}) into (\ref{mkdvPhieqs}) yields an effective characterization of $A(k)$ and $B(k)$ for the first and second Neumann problems respectively.
Indeed, substituting into the system (\ref{mkdvPhieqs}) the expansions
\begin{align*}
 \Phi_j &= \Phi_{j0} + \epsilon \Phi_{j1} + \epsilon^2 \Phi_{j2} + \cdots, \qquad j = 1,2,
  	\\
 g_j &= \epsilon g_{j1} + \epsilon^2 g_{j2} + \cdots, \qquad j = 0, 1,2,
\end{align*}
where $\epsilon > 0$ is a small parameter, we find that the terms of $O(1)$ give $\Phi_{10} \equiv 0$ and $\Phi_{20} \equiv 1$, while the terms of $O(\epsilon)$ give $\Phi_{21} \equiv 0$ and
\begin{align}\label{mkdveps1a}
&O(\epsilon): \; \Phi_{11}(t, k) = \int_0^t e^{8ik^3(t' - t)}\left(4 k^2 g_{01} + 2 i k g_{11} - g_{21}\right)dt'.
\end{align}
The symmetric combinations of equation (\ref{mkdveps1a}) are
\begin{subequations}\label{mkdvchi11}
\begin{align}\label{chi11}
& \chi_{11}(t, k) =  12 k^2 \int_0^t e^{8ik^3(t' - t)} g_{01}(t')dt',
	\\ \label{hatchi11}
& \hat{\chi}_{11}(t, k) = 6 i k  \int_0^t e^{8ik^3(t' - t)} g_{11}(t') dt',
	\\ \label{checkchi11}
& \check{\chi}_{11}(t, k) = - 3\int_0^t e^{8ik^3(t' - t)}g_{21}(t') dt'.
\end{align}
\end{subequations}

The Dirichlet problem can now be solved perturbatively as follows: Expanding (\ref{mkdvg1expression}) and (\ref{mkdvg2expression}) and assuming for simplicity that $a(k)$ has no zeros, we find
\begin{subequations}\label{mkdvg11g21expressions}
\begin{align} \nonumber
g_{11} = \;& \frac{2}{\pi i} \int_{\partial D_3}\biggl(k\chi_{11}(t,k) - \frac{3g_{01}}{2}\biggr) dk 
	\\\label{mkdvg11expression}
& - \frac{2}{\pi i} \int_{\partial D_3} \bigl[(\alpha^2 - \alpha) b_1(\alpha k) + (\alpha - \alpha^2)b_1(\alpha^2 k) \bigr] k e^{-8ik^3t} dk,
	\\   \nonumber
g_{21} = & -\frac{4}{\pi} \int_{\partial D_3}\biggl(k^2\chi_{11}(t,k) - \frac{3kg_{01}}{2}\biggr) dk 
	\\\label{mkdvg21expression}
&  + \frac{4}{\pi} \int_{\partial D_3} 
k^2 e^{-8ik^3 t} \bigl[(1-\alpha)b_1(\alpha k) + (1 - \alpha^2) b_1(\alpha^2 k)\bigr]dk.
\end{align}
\end{subequations}
Using equation (\ref{chi11}) to determine $\chi_{11}$, we can determine $g_{11}$, $g_{21}$ from (\ref{mkdvg11g21expressions}); then $\Phi_{11}$ can be found from (\ref{mkdveps1a}). These arguments, just like in the case of NLS, can be extended to higher orders and thus yields a constructive scheme for computing $A(k)$ and $B(k)$ to all orders.

Similarly, for the first Neumann problem, the expansion of (\ref{N1expressions}) gives
\begin{align}\nonumber
  g_{01}(t) =\;& \frac{1}{\pi} \int_{\partial D_3} \hat{\chi}_{11} dk 
- \frac{1}{\pi} \int_{\partial D_3} e^{-8ik^3t} \bigl[(\alpha - \alpha^2)b_1(\alpha k) + (\alpha^2 - \alpha)b_1(\alpha^2 k)\bigr]dk,
  	\\ \label{N1pertexpression}  
  g_{21}(t) =\;&- \frac{4}{\pi} \int_{\partial D_3} \left(k^2\hat{\chi}_{11} - \frac{3}{\pi i}\left(\int_{\partial D_3} l \hat{\chi}_{11}(t,l) dl\right)\right) dk 
	\\ \nonumber
& + \frac{4}{\pi} \int_{\partial D_3} 
k^2 e^{-8ik^3 t} \bigl[(1 - \alpha^2)b_1(\alpha k) + (1 - \alpha)b_1(\alpha^2 k)\bigr]dk.
\end{align}
In this case, we first use equation (\ref{hatchi11}) to determine $\hat{\chi}_{11}$; we then determine $g_{11}$, $g_{21}$ from (\ref{N1pertexpression}); then $\Phi_{11}$ can be found from (\ref{mkdveps1a}). 
For the second Neumann problem, we use a similar argument based on (\ref{checkchi11}).
In all cases, the system can be solved perturbatively to all orders.

\subsection{The mKdVII equation}
An important question in the analysis of initial-boundary value problems is the determination of the number of boundary conditions needed for a well-posed problem. For the mKdV equation (\ref{mkdvI}) considered earlier (in this subsection we will denote this equation by mKdVI), we have seen that one boundary condition must be imposed at $t = 0$.  
On the other hand, for the mKdVII equation, which has a minus sign in front of $\partial^3 q/\partial x^3$,
\begin{align}\label{mkdvII}\tag{mKdVII}
 & \frac{\partial q}{\partial t} - \frac{\partial^3 q}{\partial x^3} + 6 \lambda q^2 \frac{\partial q}{\partial x} = 0, \qquad q\;\text{real},\quad \lambda = \pm 1,
\end{align}
two boundary conditions must be prescribed at $t = 0$ for a well-posed problem. 
How does this difference between mKdVI and mKdVII manifest itself in the present formalism?

To answer this question we will assume for simplicity that $q_0 = 0$. There exist three different representations for $g_0(t)$ for mKdVI: 
\begin{align}\label{mkdvIg0reps}
g_0(t) = \frac{2}{\pi} \int_{\partial D_3} \chi_1 dk, \qquad
g_0(t) = \frac{1}{\pi} \int_{\partial D_3} \hat{\chi}_1 dk, \qquad
g_0(t) = \frac{1}{\pi} \int_{\partial D_3} \check{\chi}_1 dk.  
\end{align}  
The latter two representations are useful for the effective solutions of the first and second Neumann problems, respectively. Indeed, the perturbative equations (cf. equation (\ref{mkdvchi11}))
\begin{subequations}\label{mkdvchi1n}
\begin{align}
& \chi_{1n}(t, k) =  12 k^2 \int_0^t e^{8ik^3(t' - t)} g_{0n}(t')dt' + \text{lower order terms},
	\\
& \hat{\chi}_{1n}(t, k) = 6 i k  \int_0^t e^{8ik^3(t' - t)} g_{1n}(t') dt' + \text{lower order terms},
	\\ 
& \check{\chi}_{1n}(t, k) = - 3\int_0^t e^{8ik^3(t' - t)}g_{2n}(t') dt' + \text{lower order terms},
\end{align}
\end{subequations}
show that if $g_1$ and $g_2$ are known, then $\hat{\chi}_1$ and $\check{\chi}_1$ respectively can be determined uniquely at each step of the perturbative expansion. 

The second representation in (\ref{mkdvIg0reps}) is derived by noting that
\begin{subequations}\label{repproof}
\begin{align}\label{secondrepproof}
\int_{\partial D_3} \hat{\chi}_1 dk
= \biggl(\int_{\partial D_3} + \alpha \int_{\alpha \partial D_3} + \alpha^2 \int_{\alpha^2 \partial D_3}\biggr) \Phi_1 dk.
\end{align}
Since $\Phi_1$ is bounded and analytic in  $\C \setminus D_3$, each of the integrals on the right-hand side can be computed individually and this yields
$$\int_{\partial D_3} \hat{\chi}_1 dk = \left(\frac{5\pi i}{3} - \alpha \frac{\pi i}{3} - \alpha^2 \frac{\pi i}{3}\right) \Phi_1^{(1)}
= \pi g_0.$$
The third representation is derived in the same way by noting that
\begin{align}\label{thirdrepproof}
\int_{\partial D_3} \check{\chi}_1 dk
= \biggl(\int_{\partial D_3} + \alpha^2 \int_{\alpha \partial D_3} + \alpha \int_{\alpha^2 \partial D_3}\biggr) \Phi_1 dk.
\end{align}
\end{subequations}

Let us now consider the situation for mKdVII. We claim that the representations in (\ref{mkdvIg0reps}) have the following analogs for mKdVII:
\begin{align}\label{mkdvIIg0reps}
g_0(t) =-\frac{1}{\pi} \int_{\partial D_3} \chi_1 dk, \qquad
g_0(t) = - \frac{1}{\pi} \int_{\partial D_3}(\hat{\chi}_1 + \check{\chi}_1) dk.
\end{align}  
Indeed, mKdVII admits the Lax pair
\begin{equation}\label{mkdv2lax}
\begin{cases} \mu_x - ik[\sigma_3, \mu] = Q\mu, \\
\mu_t + 4ik^3[\sigma_3, \mu] = \tilde{Q}\mu,
\end{cases}
\end{equation}
where
$$Q(x,t) = \begin{pmatrix} 0 & q \\ \lambda q & 0 \end{pmatrix}, \qquad
\tilde{Q}(x,t, k) = -2Q^3 + Q_{xx} - 2ik(Q^2 + Q_x)\sigma_3 - 4k^2Q.$$
Proceeding as in the case of mKdVI, we find 
\begin{align*}
& \Phi_1(t, k) = -\frac{g_0(t)}{2ik} + O \Bigl(\frac{1}{k^2} \Bigr), \qquad k \to \infty, \quad k \in D_2 \cup D_4,
\end{align*}
i.e. $g_0 = -2i\Phi_1^{(1)}$. 
Using this relationship it is straightforward to derive, just like in the case of mKdVI, the first representation in (\ref{mkdvIIg0reps}). Morever, equations (\ref{repproof}) are also valid for mKdVII. However, in this case the analyticity domain $\C \setminus D_3$ of $\Phi_1$ is smaller (see figure \ref{Djsmkdv2fig}) and this prohibits the evaluation of the integrals on the right-hand sides of (\ref{repproof}). In order to proceed, we now have to consider the sum of equations (\ref{secondrepproof}) and (\ref{thirdrepproof}):
$$\int_{\partial D_3} (\hat{\chi}_1 + \check{\chi}_1) dk
= \biggl(2\int_{\partial D_3} - \int_{\alpha \partial D_3} - \int_{\alpha^2 \partial D_3}\biggr) \Phi_1 dk
= 2\pi i \Phi_1^{(1)},$$
which yields the second representation in (\ref{mkdvIIg0reps}). 

In summary, whereas there exist {\it three} different representations (\ref{mkdvIg0reps}) of $g_0$ for mKdVI, there exist only {\it two} such representations for mKdVII. Since the perturbative equations (\ref{mkdvchi1n}) are also valid for mKdVII provided that the signs of $g_{0n}$ and $g_{2n}$ are reversed, the second representation in (\ref{mkdvIIg0reps}) can be used to determine $g_0$ at each step of the perturbative scheme provided that {\it both} $g_1$ and $g_2$ are prescribed. 

Similar remarks apply to the representations of $g_1$ and $g_2$. For example, for $g_2$ there exist three different representation for mKdVI involving $\chi_1$, $\hat{\chi}_1$ and $\check{\chi}_1$, respectively; the former two are useful for the effective solutions of the Dirichlet and first Neumann problems, respectively. On the other hand, for mKdVII there exist only two representations involving $\chi_1 + \hat{\chi}_1$ and $\check{\chi}_1$ respectively (cf. equation (1.18a) in \cite{F2005}); the former can be used to give an effective solution provided that {\it both} $g_0$ and $g_1$ are prescribed.

\begin{figure}
\begin{center}
\begin{overpic}[width=.34\textwidth]{DjsmKdV.pdf}
      \put(77,57){$D_3$}
      \put(46,70){$D_4$}
      \put(15,57){$D_3$}
      \put(15,25){$D_2$}
      \put(46,12){$D_1$}
      \put(77,25){$D_2$}
       \put(101,40){$\text{Re}\; k$}
      \put(60,47){$\pi/3$}
      \end{overpic}
     \begin{figuretext}\label{Djsmkdv2fig}
       The domains $\{D_j\}_1^4$ for mKdVII.
     \end{figuretext}
     \end{center}
\end{figure}

\section{The Dirichlet to Neumann map}\label{D2Nsec}\nequation
We have shown in sections \ref{NLSsec} and \ref{mKdVIsec} that the spectral functions associated with the Dirichlet and Neumann problems for NLS and mKdV are characterized by nonlinear integral equations which can be solved perturbatively to all orders. In this section, we demonstrate that these integral equations can also be employed to construct perturbatively the generalized Dirichlet to Neumann map. 

Consider the Dirichlet problem for NLS on the half-line. The perturbative construction of the Dirichlet to Neumann map involves the following: Given initial data $q_0(x)$ and Dirichlet data $g_0(t)$ in the form
$$q_0(x) = \epsilon q_{01}(x) + \epsilon^2 q_{02}(x) + \cdots, \qquad
 g_0(t) = \epsilon g_{01}(t) + \epsilon^2 g_{02}(t) + \cdots,$$
where $\epsilon > 0$ is a small parameter, determine the coefficients $\{g_{1n}(t)\}_{n\geq1}$ in the expansion of the corresponding Neumann data,
\begin{align}\label{g1expansion}
  g_1(t) = \epsilon g_{11}(t) + \epsilon^2 g_{12}(t) + \epsilon^3 g_{13}(t) + \cdots, 
\end{align}		
in terms of $\{q_{0n}(t), g_{0n}(t)\}_{n \geq 1}$.
Assuming for simplicity that $q_0(x) \equiv 0$, we will derive expressions for $g_{11}$, $g_{12}$, and $g_{13}$ in terms of $g_{01}$, $g_{02}$, $g_{03}$ for NLS on the half-line. Similar results can be obtained for the Neumann to Dirichlet map for NLS, as well as for the various generalized Dirichlet to Neumann maps for mKdV. 

\begin{theorem}
Let 
$$q(x,t) = \epsilon q_1(x,t) + \epsilon^2 q_2(x,t) + \cdots,$$
be the perturbative solution of the NLS equation on the half-line satisfying the initial conditions $q(x,0) = 0$, $x > 0$, and the Dirichlet boundary conditions 
$$q(0,t) = \epsilon g_{01}(t) + \epsilon^2 g_{02}(t) + \cdots,$$
where $\{g_{0n}(t)\}_{n \geq 1}$ are sufficiently smooth functions compatible with the zero initial data, i.e.
$g_{0n}(0) = \dot{g}_{0n}(0) = \ddot{g}_{0n}(0) = \cdots = 0$, $n \geq 1$.
Then, the first few coefficients in the expansion (\ref{g1expansion}) of the Neumann data $g_1(t) = q_x(0,t)$ are given by the following formulas:
\begin{align} \label{g11g12formulas}
g_{11}(t) = \;&
- \frac{e^{-\frac{i\pi}{4}}}{\sqrt{\pi}} \int_0^t \frac{\dot{g}_{01}(t')}{\sqrt{t - t'}} dt',
\qquad
g_{12}(t) = \; 
- \frac{e^{-\frac{i\pi}{4}}}{\sqrt{\pi}} \int_0^t \frac{\dot{g}_{02}(t')}{\sqrt{t - t'}} dt',
	\\ \nonumber
g_{13}(t) =\; & \frac{2\lambda c}{\pi i}\int_0^t  \frac{|g_{01}(t')|^2g_{01}(t')}{\sqrt{t - t'}} dt' 
-  \frac{2\lambda c}{\pi i}\int_0^t |g_{01}(t')|^2 \int_0^{t'}\frac{\dot{g}_{01}(t'')}{\sqrt{t - t''}}dt'' dt'
	\\\nonumber
&- \frac{2c}{\pi}\int_0^t \frac{\dot{g}_{03}(t')}{\sqrt{t - t'}}dt' 
	\\\nonumber
& - \frac{\lambda c}{\pi i} \int_0^t g_{01}(t') \int_0^{t'} \bar{g}_{01}(t'') \int_0^{t''} \frac{\dot{g}_{01}(t''')}{(t-t'+t'' - t''')^{3/2}}  dt''' dt''  dt'
	\\\nonumber
&+ \frac{\lambda c}{\pi} \int_0^t g_{01}(t') \int_0^{t'} \bar{g}_{11}(t'') \int_0^{t''} \frac{g_{11}(t''')}{(t-t'+t'' - t''')^{3/2}} dt'''dt''dt'
	\\\nonumber
& - \frac{\lambda c}{\pi}\int_0^tg_{11}(t') \int_0^{t'}\bar{g}_{01}(t'') \int_0^{t''} \frac{g_{11}(t''')}{(t - t' + t'' - t''')^{3/2}}  dt'''dt'' dt'
	\\\nonumber
& +\frac{\lambda c}{\pi} \int_0^tg_{11}(t') \int_0^{t'} \bar{g}_{11}(t'')\int_0^{t''} \frac{g_{01}(t''')}{(t - t' + t'' - t''')^{3/2}}  dt'''dt''dt'
	\\ \label{g13formula}
& +   i \lambda g_{01}(t)  \int_0^t \left[\bar{g}_{01}(t')g_{11}(t') - \bar{g}_{11}(t') g_{01}(t')\right] dt',
\end{align}	
where $c = \frac{\sqrt{\pi}}{2}e^{-\frac{i\pi}{4}}$.
\end{theorem}
\proofbegin
Let
$$\Phi_1 = \epsilon \Phi_{11} + \epsilon^2 \Phi_{12} + \cdots, \qquad
\Phi_2 = \epsilon \Phi_{21} + \epsilon^2 \Phi_{22} + \cdots,$$
and let $\{\chi_j\}_1^2$ and $\{\check{\chi}_j\}_1^2$ be the odd and even combinations formed from $\{\Phi_j\}_1^2$ as in (\ref{chijdef}) and (\ref{checkchijdef}). Let
$$\chi_j = \epsilon \chi_{j1} + \epsilon^2 \chi_{j2} + \cdots, \qquad
\check{\chi}_j = \epsilon \check{\chi}_{j1} + \epsilon^2 \check{\chi}_{j2} + \cdots.$$

We first prove (\ref{g11g12formulas}). By (\ref{g11expression}),
\begin{align} \label{tperiodicg11}
g_{11}(t) = \frac{2}{\pi i}\int_{\partial D_3}\bigl(k\chi_{11}(t,k) + i g_{01}(t)\bigr) dk.
\end{align}
Substituting the expression (\ref{chi11eqn}) for $\chi_{11}$ into (\ref{tperiodicg11}) and integrating by parts, we find
\begin{align} \label{g11doubleintegral}
g_{11}(t) = \frac{2}{\pi}\int_{\partial D_3} \int_0^t e^{4ik^2(t' - t)} \dot{g}_{01}(t') dt' dk. 
\end{align}
Using the identity
\begin{align}\label{performdkidentity}
\int_{\partial D_3} e^{-4ik^2 t}dk = - \frac{c}{\sqrt{t}}, \qquad t > 0, \quad c = \frac{\sqrt{\pi}}{2}e^{-\frac{i\pi}{4}},
\end{align}
we can compute the $k$ integral in (\ref{g11doubleintegral}). This gives the first equation in (\ref{g11g12formulas}). The derivation of the second equation in (\ref{g11g12formulas}) is similar.

In order to prove (\ref{g13formula}), we note that (\ref{g1expression}) implies that
\begin{align}\label{g13twointegrals}
g_{13}(t) = \frac{2}{\pi i}\int_{\partial D_3} (k \chi_{13}(t,k) + ig_{03}(t)) dk + \frac{2g_{01}(t)}{\pi} \int_{\partial D_3} \chi_{22}(t,k) dk.
\end{align}
Moreover, by (\ref{Phieqs}),
\begin{align*} 
& \Phi_{13} = \int_0^t e^{4ik^2(t' - t)}\left[-i\lambda |g_{01}|^2\Phi_{11} + 2kg_{03} + ig_{13} + (2kg_{01} + ig_{11})\Phi_{22}\right] dt',
	\\
& \Phi_{22} = \lambda \int_0^t \left[(2k\bar{g}_{01} - i\bar{g}_{11})\Phi_{11} + i|g_{01}|^2\right] dt',
\end{align*}
so that
\begin{align*}
& \chi_{13} = \int_0^t e^{4ik^2(t' - t)}\left[-i\lambda |g_{01}|^2\chi_{11} + 4kg_{03} + 2kg_{01}\check{\chi}_{22} + ig_{11}\chi_{22}\right] dt',
	\\
& \chi_{22} = \lambda \int_0^t \left[2k\bar{g}_{01}\check{\chi}_{11} - i\bar{g}_{11}\chi_{11}\right] dt',
	\\
& \check{\chi}_{22} = \lambda \int_0^t \left[2k\bar{g}_{01}\chi_{11} - i\bar{g}_{11}\check{\chi}_{11} + 2i|g_{01}|^2\right] dt',
	\\
&  \check{\chi}_{11} = 2i \int_0^t e^{4ik^2(t' - t)}g_{11} dt'.
\end{align*}
Thus, letting $I$ denote the first integral on the right-hand side of (\ref{g13twointegrals}), i.e.
$$I = \int_{\partial D_3} (k \chi_{13}(t,k) + ig_{03}(t)) dk,$$
we find
\begin{align*}
I = &  \int_{\partial D_3} \Biggl\{k\int_0^t e^{4ik^2(t' - t)}\Biggl[-i\lambda |g_{01}(t')|^2
4k\int_0^{t'}e^{4ik^2(t'' -t')}g_{01}(t'')dt''
	\\
& + 4kg_{03}(t')
+ 2kg_{01}(t')
\lambda \int_0^{t'} \bigl[2k\bar{g}_{01}\chi_{11} - i\bar{g}_{11}\check{\chi}_{11} + 2i|g_{01}|^2\bigr](t'',k) dt''
	\\
& + ig_{11}(t')
\lambda \int_0^{t'} \bigl[2k\bar{g}_{01}\check{\chi}_{11} - i\bar{g}_{11}\chi_{11}\bigr](t'',k) dt''\Biggr] dt' + ig_{03}(t)\Biggr\} dk
	\\
 = & -i\lambda\int_{\partial D_3} \int_0^t |g_{01}(t')|^2
\int_0^{t'}4k^2e^{4ik^2(t'' -t)}g_{01}(t'')dt'' dt'dk
	\\
& + \int_{\partial D_3} \left( \int_0^t e^{4ik^2(t' - t)}4k^2g_{03}(t') dt' + ig_{03}(t)\right)dk
	\\
& + \int_{\partial D_3}  \int_0^t e^{4ik^2(t' - t)}2k^2g_{01}(t')
\lambda \int_0^{t'} \biggl[\bar{g}_{01}(t'')8k^2\int_0^{t''}e^{4ik^2(t''' -t'')}g_{01}(t''')dt''' 
	\\
&\hspace{3cm} - i\bar{g}_{11}(t'') 2i \int_0^{t''} e^{4ik^2(t''' - t'')}g_{11}(t''') dt'''
 + 2i|g_{01}(t'')|^2\biggr] dt'' dt' dk
	\\
& + \int_{\partial D_3}  \int_0^t ke^{4ik^2(t' - t)} ig_{11}(t') \lambda \int_0^{t'} \biggl[2k\bar{g}_{01}(t'')
2i \int_0^{t''} e^{4ik^2(t''' - t'')}g_{11}(t''') dt'''
	\\
& \hspace{5cm}- i\bar{g}_{11}(t'')4k\int_0^{t''}e^{4ik^2(t''' -t'')}g_{01}(t''')dt'''\biggr] dt'' dt' dk.
\end{align*}	
Integration by parts yields 
\begin{align*}
I = & -i\lambda\int_{\partial D_3} \int_0^t |g_{01}(t')|^2
\biggl(-i e^{4ik^2(t' -t)}g_{01}(t')
+i \int_0^{t'}e^{4ik^2(t'' -t)}\dot{g}_{01}(t'')dt''\biggr) dt'dk
	\\
& + \int_{\partial D_3} i\int_0^t e^{4ik^2(t' - t)}\dot{g}_{03}(t') dt' dk
	\\
& + \int_{\partial D_3}  \Biggl\{
-ig_{01}(t) \lambda \int_0^{t} \Biggl[\bar{g}_{01}(t') i\int_0^{t'}e^{4ik^2(t'' -t')}\dot{g}_{01}(t'')dt'' 
	\\
& \hspace{4.3cm} +\bar{g}_{11}(t') \int_0^{t'} e^{4ik^2(t'' - t')}g_{11}(t'') dt''\Biggr] dt'
 	\\
& \hspace{1.5cm} + i\int_0^t e^{4ik^2(t' - t)}
\frac{\partial}{\partial t'}\biggl\{g_{01}(t')\lambda \int_0^{t'} \biggl[\bar{g}_{01}(t'') i\int_0^{t''}e^{4ik^2(t''' -t'')}\dot{g}_{01}(t''')dt''' 
	\\
& \hspace{5.2cm}+\bar{g}_{11}(t'') \int_0^{t''} e^{4ik^2(t''' - t'')}g_{11}(t''') dt'''\biggr] dt'' \biggr\} dt' \Biggr\}
 dk
	\\
& + \int_{\partial D_3}  \int_0^t e^{4ik^2(t' - t)} ig_{11}(t') \lambda \int_0^{t'} \biggl[\bar{g}_{01}(t'')
\biggl( g_{11}(t'') - \int_0^{t''} e^{4ik^2(t''' - t'')}\dot{g}_{11}(t''') dt'''\biggr)
	\\
&\hspace{4cm} - \bar{g}_{11}(t'')\biggl(g_{01}(t'') - \int_0^{t''}e^{4ik^2(t''' -t'')}\dot{g}_{01}(t''')dt'''\biggr)
\biggr] dt'' dt' dk.
\end{align*}	
Using the identity (\ref{performdkidentity}) to perform the $dk$ integrals, we deduce that 
\begin{align*}
I= &\; \lambda c\int_0^t |g_{01}|^2 \biggl(\frac{g_{01}(t')}{\sqrt{t - t'}} - \int_0^{t'}\frac{\dot{g}_{01}(t'')}{\sqrt{t - t''}}dt'' \biggr) dt'
 - i c \int_0^t \frac{\dot{g}_{03}(t')}{\sqrt{t - t'}} dt' 
	\\
& + i \lambda c g_{01}(t) \int_0^{t} \biggl[ i\bar{g}_{01}(t') \int_0^{t'}\frac{\dot{g}_{01}(t'')}{\sqrt{t' - t''}}dt'' 
 +\bar{g}_{11}(t') \int_0^{t'} \frac{g_{11}(t'')}{\sqrt{t' - t''}} dt''\biggr] dt'
 	\\
&- i \lambda c \int_0^t \dot{g}_{01}(t')  \int_0^{t'} \biggl[ i\bar{g}_{01}(t'') \int_0^{t''}\frac{\dot{g}_{01}(t''')}{\sqrt{t - t' + t'' - t'''}}dt''' 
	\\
& \hspace{4cm} +\bar{g}_{11}(t'') \int_0^{t''} \frac{g_{11}(t''')}{\sqrt{t - t' + t'' - t'''}} dt'''\biggr] dt'' dt' 
 	\\
&- i \lambda c \int_0^t g_{01}(t')  \biggl[ i\bar{g}_{01}(t')\int_0^{t'} \frac{\dot{g}_{01}(t'')}{\sqrt{t- t''}} dt'' 
 +\bar{g}_{11}(t') \int_0^{t'} \frac{g_{11}(t'')}{\sqrt{t- t''}} dt''
 \biggr] dt' 
	\\
& - i \lambda c \int_0^t \frac{g_{11}(t')}{\sqrt{t - t'}} \int_0^{t'} \bigl[\bar{g}_{01}(t'')g_{11}(t'') - \bar{g}_{11}(t'')g_{01}(t'')\bigr]dt'' dt'
	\\
& + i\lambda c \int_0^tg_{11}(t') \int_0^{t'}\bar{g}_{01}(t'') \int_0^{t''} \frac{\dot{g}_{11}(t''')}{\sqrt{t - t' + t'' - t'''}} dt'''dt'' dt'
	\\
& - i \lambda c\int_0^tg_{11}(t') \int_0^{t'} \bar{g}_{11}(t'')\int_0^{t''} \frac{ \dot{g}_{01}(t''')}{\sqrt{t - t' + t'' - t'''}} dt'''dt''dt'.
\end{align*}	
We next integrate the terms involving $1/\sqrt{t - t' + t'' - t'''}$ by parts with respect to $dt'$. The total contribution of lines $2$-$5$ of the above expression yields
\begin{align*}
&i\lambda c \int_0^t g_{01}(t') \int_0^{t'} \biggl[i \bar{g}_{01}(t'')\int_0^{t''} \frac{\dot{g}_{01}(t''')}{2(t-t'+t'' - t''')^{3/2}} dt'''
	\\
&\hspace{3.2cm} + \bar{g}_{11}(t'') \int_0^{t''} \frac{g_{11}(t''')}{2(t-t'+t'' - t''')^{3/2}} dt'''\biggr]dt''dt',
\end{align*}
while the last two lines combine nicely with the sixth line. We thus arrive at
\begin{align}\nonumber
I = &\; \lambda c \int_0^t |g_{01}|^2 \left(\frac{g_{01}(t')}{\sqrt{t - t'}} - \int_0^{t'}\frac{\dot{g}_{01}(t'')}{\sqrt{t - t''}}dt'' \right) dt' 
- i c \int_0^t \frac{\dot{g}_{03}(t')}{\sqrt{t - t'}} dt' 
	\\\nonumber
& -\lambda c \int_0^t g_{01}(t') \int_0^{t'} \bar{g}_{01}(t'') \int_0^{t''} \frac{\dot{g}_{01}(t''')}{2(t-t'+t'' - t''')^{3/2}}  dt'''dt''dt'
	\\\nonumber
&+  i \lambda c \int_0^t g_{01}(t') \int_0^{t'} \bar{g}_{11}(t'') \int_0^{t''} \frac{g_{11}(t''')}{2(t-t'+t'' - t''')^{3/2}} dt'''dt''dt'
	\\\nonumber
& - \lambda c \int_0^tig_{11}(t') \int_0^{t'}\bar{g}_{01}(t'') \int_0^{t''} \frac{g_{11}(t''')}{2(t - t' + t'' - t''')^{3/2}}  dt'''dt'' dt'
	\\\label{firstintegralexpression}
& + \lambda c \int_0^tig_{11}(t') \int_0^{t'} \bar{g}_{11}(t'') \int_0^{t''} \frac{g_{01}(t''')}{2(t - t' + t'' - t''')^{3/2}}  dt'''dt''dt'.
\end{align}	

On the other hand, using the identity
$$\int_{\partial D_3} k \int_0^t e^{4ik^2(t'' - t')}f(t'')dt'' dk = \begin{cases} \frac{\pi}{4}f(t'), & 0 < t' < t, \\
\frac{\pi}{8}f(t'), & 0 < t' = t, 
\end{cases}$$
where $f(t'')$ is a smooth function, we find that the second integral on the right-hand side of (\ref{g13twointegrals}) can be written as follows:
\begin{align}\nonumber
 \int_{\partial D_3} \chi_{22}(t,k) dk
& = \int_{\partial D_3} \lambda \int_0^t \biggl[2k\bar{g}_{01}(t') 2i \int_0^{t'} e^{4ik^2(t'' - t')}g_{11}(t'') dt'' 
	\\ \nonumber
& \hspace{3cm} - i\bar{g}_{11}(t')4k\int_0^{t'} e^{4ik^2(t'' -t')}g_{01}(t'')dt''\biggr] dt'dk
	\\ \label{secondintegralexpression}
& = i \lambda \frac{\pi}{2} \int_0^t \left[\bar{g}_{01}(t')g_{11}(t') 
- \bar{g}_{11}(t') g_{01}(t')\right] dt'.
\end{align}
Equations (\ref{g13twointegrals}), (\ref{firstintegralexpression}), and (\ref{secondintegralexpression}) imply (\ref{g13formula}).
\proofend

\section{Conclusions}\label{conclusionsec}\nequation
As it was mentioned in the introduction, linear PDEs involving only second order derivatives can be analyzed by an appropriate extension from the half-line to the full line. In analogy with the linear case, linearizable BVPs for the NLS and the sG equations have been studied via techniques based on an appropriate extension from the half-line to the real line. These extensions yield explicit conditions on the scattering data for an initial value problem formulated on the full line, see \cite{AS1975, BT1991, F1989, H1990, T1988}.

It is of course possible to map the RH problem obtained by the unified method to the RH problem obtained via an extension from the half-line to the full line, see for example \cite{FK2004}. However, as it was noted in the introduction, the unified method has the advantage that it always yields solutions that are uniformly convergent at the boundaries. Furthermore, the unified method is apparently the {\it only} approach that can be applied to linearizable BVPs for PDEs involving \textit{third} order derivatives such as the KdV and the mKdV equations.

\appendix
\section{A comparison of the unified method with the \\ approach of \cite{DMS2001}} \label{DMSapp}
\renewcommand{\theequation}{A.\arabic{equation}}\nequation
In sections \ref{NLSsec} and \ref{mKdVIsec} an effective characterization of the generalized Dirichlet to Neumann map for the NLS and the mKdV equations respectively, was presented. The effectiveness was demonstrated by showing that the relevant nonlinear equations can be solved perturbatively via a well-defined recursive scheme. Furthermore, in the linear limit, the relevant formulas coincide with the formulas obtained by solving, via the new method, the associated linearized equation. 

It was noted in the introduction that for the special case of the NLS on the half-line (but {\it not} for the mKdV on the half-line, or even for the NLS on the finite interval), the approach introduced in \cite{DMS2001} yields expressions for $g_0(t)$ in terms of $q_0(x)$ and $g_1(t)$ for the Neumann problem and similarly for the Dirichlet and the Robin problems.
In what follows we compare the unified method with the approach of \cite{DMS2001} for the following slight generalization of the linearized version of the NLS:
\begin{align}\label{appnls}
  iu_t(x,t) + u_{xx}(x,t) + i\alpha u_x(x,t) = 0, \qquad 0 < x < \infty, \; 0 < t < T,
\end{align}
where $\alpha$ is a real constant. 

\subsection{The unified method}
We first implement the new method. Equation (\ref{appnls}) can be written in the following divergence form:
\begin{subequations}
\begin{align}\label{appnlsdivform}
  \left(ue^{-ikx + i w(k)t} \right)_t - \left[(iu_x - (\alpha + k)u)e^{-ikx + iw(k)t}\right]_x = 0, \qquad k \in \C,
  \end{align}
  where the dispersion relation $w(k)$ is given by
 \begin{align}
 w(k) = k^2 + \alpha k.
 \end{align}
 \end{subequations}
 We note that equation (\ref{appnlsdivform}) and Green's theorem immediately imply the global relation
 \begin{align}\label{GRappA}
 e^{iw(k)t}\hat{u}(k,t) = \hat{u}_0(k) - \tilde{g}(k,t), \qquad \text{Im}\, k \leq 0,
 \end{align}  
where $\hat{u}(k,t)$ denotes the half-Fourier transform of $u(x,t)$ and the spectral functions $\{\hat{u}_0(k), \tilde{g}(k,t)\}$ are defined as follows:
\begin{align}
\hat{u}_0(k) = \int_0^\infty u_0(x) e^{-ikx} dx, \qquad k \in \C^-,
\end{align}
\begin{subequations}
  \begin{align}
  &  \tilde{g}(k,t) = i\tilde{g}_1(w(k),t) - (\alpha + k)\tilde{g}_0(w(k),t), \qquad k \in \C,
    	\\
   & \tilde{g}_j(k,t) = \int_0^t e^{ik\tau} \partial_x^ju(0,\tau)d\tau, \qquad j = 0,1.
  \end{align}
\end{subequations}
We note that equation (\ref{GRappA}) can also be obtained by the usual application of the half-Fourier transform.

Employing the inverse Fourier transform, equation (\ref{GRappA}) yields
\begin{align}\label{A6}
  u(x,t) = \frac{1}{2\pi} \int_{-\infty}^\infty e^{-ikx - iw(k)t} \hat{u}_0(k) dk
  - \frac{1}{2\pi} \int_{-\infty}^\infty e^{-ikx - iw(k)t} \tilde{g}(k,t) dk.
\end{align}
Jordan's lemma implies that it is possible to replace in the second term of (\ref{A6}) the integral along the real axis with an integral along the contour $\partial D$ which consists of the union of two rays:
$$\partial D: \quad \biggl\{\text{Im}\, k = 0, -\frac{\alpha}{2} < \re k < \infty\biggr\} \cup \biggl\{ \re k = -\frac{\alpha}{2}, 0 < \text{Im}\, k < \infty\biggr\}.$$

We next consider the transformations $k \to \nu(k)$, which leave $w(k)$ invariant:
$$k^2 + \alpha k = \nu^2 + \alpha \nu, \quad \text{i.e.}\quad \nu = k \; \text{and} \; \nu = -k - \alpha.$$
Replacing in equation (\ref{GRappA}) $k$ with $\nu(k) = -k -\alpha$, we find
\begin{align}\label{A7}
  e^{iw(k)t} \hat{u}(-k-\alpha, t) = \hat{u}_0(-k - \alpha) - i \tilde{g}_1 - k\tilde{g}_0 = 0, \qquad \text{Im}\, \geq 0.
 \end{align}
 In the case of the Neumann problem, we use (\ref{A7}) to eliminate $\tilde{g}_0$, i.e. to eliminate the {\it transform} of the unknown Dirichlet boundary values. Solving equation (\ref{A7}) for $\tilde{g}_0$, substituting the resulting equation in (\ref{A6}), where the contour in the second integral is given by $\partial D$, and employing Jordan's lemma to show that the contribution of $\hat{u}(-k-\alpha, t)$ vanishes, we find
 \begin{align} \label{A8}
   u(x,t) = &\; \frac{1}{2\pi} \int_{-\infty}^\infty e^{-ikx - iw(k)t}\hat{u}_0(k)dk
   	\\\nonumber
&   - \frac{1}{2\pi} \int_{\partial D} e^{-ikx - iw(k)t} \left[i\biggl(2 + \frac{1}{k}\biggr)\tilde{g}_1(w(k),t) - \biggl(1 + \frac{\alpha}{k}\biggr)\hat{u}_0(-k-t)\right]dk.
 \end{align}

\begin{remark}\upshape
  1. Equation (\ref{A8}) expresses $u(x,t)$ in terms of the Fourier transform $\hat{u}_0(k)$ of the initial data $u_0(x)$ and of a $t$-transform of the Neumann data $g_1(t) := \partial_xu(0,t)$.
   
  2. By employing Jordan's lemma, it follows that $\tilde{g}_1(w(k), t)$ can be replaced with $\tilde{g}_1(w(k))$, where
 \begin{align}
 \tilde{g}_1(w(k)) = \int_0^T e^{iw(k)\tau} g_1(\tau)d\tau, \qquad k \in \C.
 \end{align}
\end{remark}

\subsection{The approach of \cite{DMS2001}}
We next consider the approach of \cite{DMS2001} in the special case of $\alpha = 0$. In this case, equation (\ref{A6}) becomes
\begin{align}\label{A10}
u(x,t) = &\; \frac{1}{2\pi} \int_{-\infty}^\infty e^{-ikx - iw(k)t} \hat{u}_0(k) dk 
	\\
& - \frac{1}{2\pi} \int_{-\infty}^\infty e^{-ikx - iw(k)t} \left[i\tilde{g}_1(w(k), t) - k\tilde{g}_0(w(k), t)\right] dk.
\end{align}
The important new idea introduced in \cite{DMS2001} has the following implementation for linear PDEs: Eliminate {\it directly} $\tilde{g}_0$ from equation (\ref{A10}). In this respect, evaluating equation (\ref{A10}) at $x = 0$, multiplying the resulting equation by $\exp(ik^2\tau)$, and integrating with respect to $\tau$ from $\tau = 0$ to $\tau = t$, we find
\begin{align}\label{A11}
  \tilde{g}_0(k^2, t) = \frac{1}{2\pi} \int_0^t e^{ik^2 \tau} \left\{ \int_{-\infty}^\infty e^{-il^2 \tau}\hat{u}_0(l) dl
   - \int_{-\infty}^\infty e^{-il^2\tau} i \tilde{g}_1(l^2, \tau)dl - I_0(\tau)\right\} d\tau,
\end{align} 
where
\begin{align}
  I_0(\tau) = \int_{-\infty}^\infty l e^{-il^2 \tau}\tilde{g}_0(l^2, \tau) dl.
\end{align} 
Due to the invariance of $\tilde{g}_0(l^2,\tau)$ under the transformation $l \to -l$, $I_0(\tau)$ vanishes and (\ref{A11}) becomes
\begin{align}\label{A13}
  \tilde{g}_0(k^2, t) = \frac{1}{2\pi} \int_0^t e^{ik^2 \tau} \left\{ \int_{-\infty}^\infty e^{-il^2 \tau}\hat{u}_0(l) dl
   - \int_{-\infty}^\infty e^{-il^2\tau} i \tilde{g}_1(l^2, \tau)dl\right\} d\tau.
\end{align} 
This equation expresses directly $\tilde{g}_0$ in terms of the known transforms $\hat{u}_0(l)$ and $\tilde{g}_1(l^2, \tau)$. Substituting (\ref{A13}) into (\ref{A10}) we find a formula expressing $u(x,t)$ in terms of the given initial and boundary conditions (albeit, this formula is more complicated than equation (\ref{A8}) evaluated at $\alpha = 0$).

It appears that the above approach cannot be applied to PDEs whose dispersion relation is not invariant under the transformation $k \to -k$. For example, if $\alpha \neq 0$, the analogue of equation (\ref{A11}) is now the following equation:
\begin{align} \nonumber
  \tilde{g}_0(w(k), t) = \frac{1}{2\pi} \int_0^t e^{iw(k) \tau} \biggl\{ &\int_{-\infty}^\infty e^{-iw(l) \tau}\hat{u}_0(l) dl
 	\\ \label{A14}
&   - \int_{-\infty}^\infty e^{-iw(l)\tau} i \tilde{g}_1(w(l), \tau)dl - I(\tau)\biggr\} d\tau,
\end{align} 
where
$$I(\tau) = \int_{-\infty}^\infty e^{-iw(l)\tau}(\alpha + l) \tilde{g}_0(w(l), \tau) dl.$$
In this case $I(\tau) \neq 0$, hence (\ref{A14}) does {\it not} provide an explicit expression for $\tilde{g}_0$ in terms of the given initial and boundary conditions. 

The above discussion suggests that the effective treatment of linear and integrable nonlinear PDEs requires the explicit use of the transformations which leave the linearized dispersion relation invariant.

\section{The asymptotics of $c(t,k)$} \label{clemmaapp}
\renewcommand{\theequation}{B.\arabic{equation}}\nequation
We will prove lemma \ref{clemma}; the proof of lemmaÊ \ref{mKdVclemma} is similar.
The functions $\{\Phi_j\}_1^2$ satisfy
\begin{align}\label{hatPhi12system}
 & \Phi_{1t} = -4ik^2\Phi_1 - i\lambda |g_0|^2 \Phi_1 + (2kg_0 + ig_1)\Phi_2,
  	\\ \nonumber
 & \Phi_{2t} = \lambda(2k\bar{g}_0 - i\bar{g}_1) \Phi_1 + \lambda i |g_0|^2 \Phi_2.
\end{align}
It follows that these functions admit an expansion of the form (see chapter 6 of \cite{CL1955})
\begin{align}\label{Phihatalphabeta}
\begin{pmatrix} \Phi_1(t,k) \\ \Phi_2(t,k) \end{pmatrix}
&= \Bigl(\alpha_0(t) + \frac{\alpha_1(t)}{k} + \cdots\Bigr)
+ \Bigl(\beta_0(t) + \frac{\beta_1(t)}{k} + \cdots\Bigr)e^{-4ik^2t}, 
	\\ \nonumber
& \hspace{7cm} k \to \infty, \; k \in \C,
\end{align}
where the coefficients $\alpha_j(t)$, $\beta_j(t)$, $j \geq 0$, are column vectors which are independent of $k$. 
We determine the coefficients by substituting (\ref{Phihatalphabeta}) into (\ref{hatPhi12system})
and using the initial conditions
$$\alpha_0(0) + \beta_0(0) = (0,1)^T, \qquad
\alpha_1(0) + \beta_1(0) = (0,0)^T.$$
This yields
\begin{align}\nonumber
\begin{pmatrix} \Phi_1(t,k) \\ \Phi_2(t,k) \end{pmatrix}
= & \begin{pmatrix} 0 \\ 1 \end{pmatrix} + \begin{pmatrix} \Phi_1^{(1)}(t) \\ \Phi_2^{(1)}(t) \end{pmatrix} \frac{1}{k} 
+ \begin{pmatrix} \Phi_1^{(2)}(t) \\ \Phi_2^{(2)}(t) \end{pmatrix} \frac{1}{k^2} 
+ O\Bigl(\frac{1}{k^2}\Bigr)
	\\\label{Phi1Phi2expansion}
& + \left[-\begin{pmatrix} \Phi_1^{(1)}(0) \\ 0 \end{pmatrix}\frac{1}{k} + O\Bigl(\frac{1}{k^2}\Bigr)\right]e^{-4ik^2t},
\qquad k \to \infty, \; k \in \C.
\end{align}
Substituting these expansions into the global relation (\ref{GRc}), we find 
\begin{align}\label{calmostthere}
& c(t, k) = \frac{\Phi_1^{(1)}(t)}{k} + \frac{\Phi_1^{(2)}(t)}{k^2} + O \Bigl(\frac{1}{k^3} \Bigr) 
 + \biggl(\frac{b^{(1)} - \Phi_1^{(1)}(0)}{k} + O \Bigl(\frac{1}{k^2} \Bigr) \biggr)e^{-4ik^2t} , 
 	\\ \nonumber
&\hspace{9cm}  k \to \infty, \quad k \in \C^+.
\end{align}
The assumption that $c(t,k)$ is of $O(1/k)$ as $k \to \infty$ in $\C^+$ implies that the terms in (\ref{calmostthere}) involving $e^{-4ik^2t}$ must vanish, i.e. for consistency we require
\begin{align}\label{b1Phi11}
b^{(1)} = \Phi_1^{(1)}(0).
\end{align}
Using the expressions for $b^{(1)}$ and $ \Phi_1^{(1)}(0)$, we find that equation (\ref{b1Phi11}) is valid iff the initial and boundary conditions are compatible, i.e. iff
\begin{align*}
q_0(0) = g_0(0).
\end{align*}

\begin{example}\upshape
When $\lambda = -1$, NLS admits the one-soliton solution
$$q(x,t) = e^{i x} \text{sech}(x - 2 t).$$
For this solution it is possible to compute $\{\Phi_j\}_1^2$ explicitly:
\begin{align*}
 & \Phi_1(t,k) = \frac{ e^{-4 i k^2 t} \left(i (2 k+1) \left(-2 e^{2 t+4 i k^2 t}+e^{4 t}+1\right)-e^{4 t}+1\right)}{\left((2 k+1)^2+1\right) \left(e^{4 t}+1\right)},
  	\\
 & \Phi_2(t,k) = \frac{\frac{2 e^{2 t-4 i k^2 t}}{e^{4 t}+1}+(2 k+1) (2 k-i \tanh (2 t)+1)}{(2 k+1)^2+1}.
\end{align*}
It is easy to verify that the relevant asymptotics is of the form (\ref{Phi1Phi2expansion}).
\end{example}

\section{The necessity of computing $A(k)$ and $B(k)$} \label{ABapp}
\renewcommand{\theequation}{C.\arabic{equation}}\nequation
We consider NLS with vanishing initial conditions and `small' Dirichlet boundary conditions,
\begin{align}\label{}
q(0,t) = \epsilon g_{01}(t) + \epsilon^2 g_{02}(t) + \epsilon^3 g_{03}(t) + O(\epsilon^4), \qquad \epsilon \to 0.
\end{align}
We expand the unknown Neumann boundary values in the form
\begin{align}
  q_x(0,t) = \epsilon g_{11}(t) + \epsilon^2 g_{12}(t) + \epsilon^3 g_{13}(t) + O(\epsilon^4), \qquad \epsilon \to 0.
\end{align}
The global relation (\ref{globalrelation}) becomes
\begin{align}\label{ABglobalrelation}
  A(k) = \overline{F_2(T, \bar{k})}, \qquad B(k) = -e^{4ik^2T}\overline{F_1(T,\bar{k})}, \qquad k \in \C^+,
\end{align}
where $F_1$ and $F_2$ are analytic functions of $k$ in $\C^+$ and of order $1 + O(1/k)$ and $O(1/k)$ respectively as $k \to \infty$ in $\C^+$.

The spectral functions can be expanded in the form
\begin{align} \nonumber
  A(k) &= 1 + \epsilon^2 A_2(k) + O(\epsilon^3), 
  	\\
  B(k)& = \epsilon B_1(k) + \epsilon^2 B_2(k) + \epsilon^3 B_3(k) + O(\epsilon^4), \qquad \epsilon \to 0, \; k \in \C.
\end{align}
The definition of $B(k)$ yields
\begin{align}
  B_1(k) = i\int_0^T e^{4ik^2\tau}g_{11}(\tau)d\tau - 2k\int_0^T e^{4ik^2\tau} g_{01}(\tau)d\tau, \qquad k \in \C.
\end{align}
This implies the important symmetry relation
\begin{align}\label{B1symmetry}
  B_1(-k) = B_1(k) + 4k\int_0^T e^{4ik^2\tau}g_{01}(\tau)d\tau, \qquad k \in \C.
\end{align}
Replacing in the $O(\epsilon)$ terms of the global relation (\ref{ABglobalrelation}b) $k$ with $-k$ and then using the symmetry relation (\ref{B1symmetry}), we find
\begin{align}
  B_1(k) = -4k\int_0^T e^{4ik^2\tau} g_{01}(\tau)d\tau - e^{4ik^2T} F_1(T, -k), \qquad k \in \C^-.
\end{align}
The function $F_1(T, -k)$ is unknown, however using the fact that it is bounded and analytic in $k$ for Ê$k \in \C^-$, it follows that its contribution vanishes.

It turns out that the situation for $B_2(k)$ is similar, however, the unknown functions entering in the right-hand side of the equation for $B_3(k)$ yield a {\it nontrivial contribution}. Indeed, in order to express $\{A(-k), B(-k)\}$ in terms of $\{A(k), B(k)\}$ and $g_0(t)$, we replace $k$ by $-k$ in the ODE (\ref{Phieq}) satisfied by $\Phi$; denoting by  $\; \hat{}\; $  the operation of replacing $k$Ê with  $-k$ we find:
\begin{align}\label{Phieqhat}
  \hat{\Phi}_t + 2ik^2\hat{\sigma}_3 \hat{\Phi} = (-2k\tilde{Q}_0 + \tilde{Q}_0^{(1)})\hat{\Phi}.
\end{align}
Hence, equations (\ref{Phieq}) and (\ref{Phieqhat}) imply 
\begin{align}\label{hatPhiPhieq}
  \left(e^{2ik^2t\hat{\sigma}_3} \hat{\Phi}^{-1}\Phi\right)_t = 4ke^{2ik^2t \hat{\sigma}_3} \hat{\Phi}^{-1}\tilde{Q}_0\Phi, \qquad k \in \C.
\end{align}
The $(12)$ component of this equation yields
\begin{align}
  \left(\hat{B}\bar{A} - \hat{\bar{A}}B\right)_t = 4k\left(e^{4ik^2t}g_0\bar{A}\hat{\bar{A}} - \lambda e^{-4ik^2t}\bar{g}_0B\hat{B}\right), \qquad k \in \C.
\end{align}
The $O(\epsilon)$ term of this equation yields equation (\ref{B1symmetry}), whereas the $O(\epsilon^2)$ and $O(\epsilon^3)$ terms yield the following equations:
\begin{align}\label{B2eq}
  \hat{B}_2 = B_2 + 4k\int_0^T e^{4ik^2\tau}g_{02}(\tau)d\tau,Ê\qquad k \in \C,
\end{align}
\begin{subequations}\label{B3eq}
\begin{align}
  \hat{B}_3 = B_3 + 4k\int_0^Te^{4ik^2\tau}g_{03}(\tau)d\tau + U(k), \qquad k \in \C,	
\end{align}
where
\begin{align}
  U(k) = B_1\hat{\bar{A}}_2 - \hat{B}_1 \bar{A}_2 + 4k\int_0^T\left[e^{4ik^2\tau}g_{01}(\bar{A}_2 + \hat{\bar{A}}_2) - e^{-4ik^2\tau}\bar{g}_{01} B_1 \hat{B}_1\right] d\tau.
\end{align}
\end{subequations}
Replacing in the $O(\epsilon^2)$ and $O(\epsilon^3)$ terms of the global relation (\ref{ABglobalrelation}b) $k$ with $-k$ and using equations (\ref{B2eq}) and (\ref{B3eq}), we find
\begin{align}
 & B_2 = -4k\int_0^T e^{4ik^2\tau}g_{02}(\tau)d\tau - e^{4ik^2T}F_2(T, -k), \qquad k \in \C^-,
  	\\
 & B_3 = -4k\int_0^T e^{4ik^2\tau}g_{03}(\tau)d\tau - e^{4ik^2T}F_3(T, -k) - U(k), \qquad k \in \C^-.	
\end{align}
The contribution of the terms involving $F_2(T, -k)$ and $F_3(T, -k)$ vanishes, however, $U(k)$ involves terms, such as $B_1$, which are {\it not} bounded for $k \in \C^-$, thus the contribution from $U(k)$ does {\it not} vanish. This makes it necessary to determine $B_1$ in terms of $g_0$.

\section{A less effective approach}
\renewcommand{\theequation}{D.\arabic{equation}}\nequation
It is possible to characterize $F_1(t,k)$ and $F_2(t,k)$ in terms of a system of nonlinear singular integrodifferential equation. Indeed, the functions $\hat{F}$Ê and $F$ satisfy equation (\ref{hatPhiPhieq}) but only for $k \in \R$. The Ê$(12)$ component of this equation yields
\begin{align}\label{hatFsymmequation}
  e^{4ik^2t} \left(\hat{F}_1 F_2 - F_1\hat{F}_2\right)_t = 4ke^{4ik^2t}\left(g_0 F_1 \hat{F}_1 - \lambda \bar{g}_0 F_2 \hat{F}_2\right), \qquad k \in \R.
\end{align}
Assuming that $F_1(t,k)$ has no zeros for $k \in \C^+$, the determinant condition
$$F_2 \bar{F}_2 - \lambda F_1 \bar{F}_1 = 1, \qquad k \in \R,$$
yields
\begin{align}\label{F2fromdeterminant}
  F_2(t,k) = e^{\frac{1}{2\pi i}\int_{-\infty}^\infty \ln(1 + |F_1|^2(t,l))\frac{dl}{l -k}} , \qquad k \in \C^+.
\end{align}
The function $F_1(t,k)$ is bounded and analytic for $k \in \C^-$, hence it can be represented in the form
\begin{align}\label{F1fHilberttransforms}
  F_1 = f + iHf, \qquad \bar{F}_1 = f - iHf, \qquad \hat{F_1} = \hat{f} - iHf,
\end{align}
where
\begin{align}
  (Hf)(k) = -\frac{1}{\pi} \dashint_{-\infty}^\infty \frac{f(l)dl}{l -k}, \qquad k \in \R.
\end{align}
Substituting equations (\ref{F2fromdeterminant}) and (\ref{F1fHilberttransforms}) in equation (\ref{hatFsymmequation}) and equating the real and imaginary parts of the resulting equations, we obtain two equations for the two unknown functions $f$ and $\hat{f}$.

\begin{remark}\upshape
1. The $(11)$ component of equation (\ref{hatPhiPhieq}) yields
$$\left(\hat{F}_1\bar{F}_1 - \lambda \hat{F}_2 \bar{F}_2\right)_t = 4\lambda k\left(g_0\hat{F}_1\bar{F}_2 - \bar{g}_0 \bar{F}_1 \hat{F}_2\right), \qquad k \in \C^-.$$
Each term of this equation is bounded and analytic for $k \in \C^-$, thus the substitutions (\ref{F1fHilberttransforms}) yield a {\it single} equation, which actually turns out to be an identity.

2. The appearance of the Hilbert transform makes the above system of equations hard to analyze. A system of singular integrodifferential equations was also presented in \cite{DMS2005}, see equation (4.6). Actually the first such system of equations was presented in \cite{F1989} by using an odd extension from the half-line to the full line.
It appears that the system of \cite{F1989} is simpler than both the system presented in \cite{DMS2005} and the system presented here.
\end{remark}

\bigskip
\noindent
{\bf Acknowledgement} {\it The authors acknowledge support from the EPSRC, UK. ASF acknowledges support from the Guggenheim foundation, USA.}

\bibliographystyle{plain}
\bibliography{is}

\begin{thebibliography}{99}
\small

\bibitem{AS1975}
M. J. Ablowitz and H. Segur, The inverse scattering transform: semi-infinite interval, {\it J. Math. Phys.} {\bf 16} (1975), 1054--1056.

\bibitem{AF2005}
Y. A. Antipov and A. S. Fokas, The modified Helmholtz equation in a semi-strip, {\it Math. Proc. Cambridge Philos. Soc.} {\bf 138} (2005), 339--365. 

\bibitem{bAF1999}
C. ben-Avraham and A. S. Fokas, The solution of the modified Helmholtz equation in a wedge and an application to diffusion-limited coalescence {\it Phys. Lett. A} {\bf 263} (1999), 355--359.

\bibitem{bAF2001} 
D. ben-Avraham and A.S. Fokas, The modified Helmholtz equation in a triangular domain and an application to diffusion-limited coalescence, {\it Phys. Rev. E} {\bf 64} (2001), 016114--6.

\bibitem{BT1991}
R. F. Bikbaev and V. O. Tarasov, Initial boundary value problem for the nonlinear Schr\"odinger equation, {\t J. Phys. A} {\bf 24} (1991), 2507--2516.

  \bibitem{BF2008}
J. L. Bona, and A. S. Fokas, Initial-boundary-value problems for linear and integrable nonlinear dispersive partial differential equations, {\it Nonlinearity} {\bf 21} (2008), T195--T203.
  
\bibitem{BFS2003}
A. Boutet de Monvel, A. S. Fokas, and D. Shepelsky, The analysis of the
global relation for the nonlinear Schr\"odinger equation on the half-line,
{\it Lett. Math. Phys.} {\bf 65} (2003), 199--212.

\bibitem{BFS2004}
A. Boutet de Monvel, A. S. Fokas, and D. Shepelsky, The mKdV equation on the half-line,
{\it J. Inst. Math. Jussieu} {\bf 3} (2004), 139--164. 

\bibitem{BFS2006}
A. Boutet De Monvel, A. S. Fokas, and D. Shepelsky, Integrable nonlinear evolution equations on a finite interval,
{\it Comm. Math. Phys.} {\bf 263} (2006), 133--172.

\bibitem{BIK2007}
A. Boutet de Monvel, A. Its, and V. Kotlyarov, Long-time asymptotics for the focusing NLS equation with time-periodic boundary condition, 
{\it C. R. Math. Acad. Sci. Paris} {\bf 345} (2007), 615--620.

\bibitem{BIK2009}
A. Boutet de Monvel, A. Its, and V. Kotlyarov, Long-time asymptotics for the focusing NLS equation with time-periodic boundary condition on the half-line, {\it Comm. Math. Phys.} {\bf 290} (2009), 479--522.

\bibitem{BKSZ2010}
A. Boutet de Monvel, V. Kotlyarov, D. Shepelsky, and C. Zheng, Initial boundary value problems for integrable systems: towards the long time asymptotics, {\it Nonlinearity} {\bf 23} (2010), 2483.

\bibitem{BS2003}
A. Boutet de Monvel and D. Shepelsky, The modified KdV equation on a finite interval,
{\it C. R. Math. Acad. Sci. Paris} {\bf 337} (2003), 517--522.

\bibitem{CL1955}
E. A. Coddington and N. Levinson, {\it Theory of differential equations}, New York, McGraw-Hill, 1955.

\bibitem{CF2004} 
D. Crowdy and A. S. Fokas, Explicit integral solutions for the plane elastostatic semi-strip, 
{\it Proc. R. Soc. London A} {\bf 460} (2004), 1285--1309.
  
\bibitem{D2007}
G. Dassios, What non-linear methods offered to linear problems? The Fokas transform method, {\it Int. J. Nonl. Mech.} {\bf 42} (2007), 146--156.

\bibitem{DF2005}
G. Dassios and A. S. Fokas, The basic elliptic equations in an equilateral triangle, {\it Proc. R. Soc. Lond. Ser. A} {\bf 461} (2005), 2721--2748.
 
  \bibitem{DMS2001}
A. Degasperis, S. V. Manakov, and P. M. Santini, On the initial-boundary value problems for soliton equations, {\it JETP Letters} {\bf 74} (2001), 481--485.

  \bibitem{DMS2005}
A. Degasperis, S. V. Manakov, and P. M. Santini, Integrable and nonintegrable initial boundary value problems for soliton equations, {\it J.  Nonl. Math. Phys.} {\bf 12} (2005), 228--243.

\bibitem{DVZ1994}
P. Deift, S. Venakides, and X. Zhou, The collisionless shock region for the long time behavior of the solutions of the KdV equation, {\it Comm. Pure Appl. Math.} {\bf 47} (1994), 199--206.

\bibitem{DVZ1997}
P. Deift, S. Venakides, and X. Zhou, New results in small dispersion KdV by an extension of the steepest descent method for Riemann-Hilbert problems, {\it Int. Math. Res. Not.} {\bf 1997} (1997), 286--299.

\bibitem{DZ1992}
P. A. Deift and X. Zhou, A steepest descent method for oscillatory Riemann-Hilbert problems, 
{\it Bull. Am. Math. Soc.} {\bf 20} (1992), 119--123.

\bibitem{DZ1993}
P. Deift and X. Zhou, A steepest descent method for oscillatory Riemann-
Hilbert problems. Asymptotics for the mKdV equation, 
{\it Ann. Math.} {\bf 137} (1993), 295--368.

\bibitem{DTVpreprint}
B. Deconinck, T. Trogdon, and V. Vasan, Solving linear partial differential equations, preprint.
 
\bibitem{D2009}
G. M. Dujardin, Asymptotics of linear initial boundary value problems with periodic boundary data on the half-line and finite intervals, {\it Proc. R. Soc. Lond. A} {\bf 465} (2009), 3341--3360.

\bibitem{D2010}
M. Doschoris, Harmonic functions in rectangular domains. Classical solutions revisited, {\it Int. J. Math. Anal.} {\bf 4} (2010), 2261--2285.

\bibitem{FF2008}
N. Flyer and A. S. Fokas, A hybrid analytical numerical method for solving evolution partial differential equations. I. The half-line, {\it Proc. R. Soc. A} \textbf{464} (2008), 1823--1849.

\bibitem{F1989}
A. S. Fokas, {An initial-boundary value problem for the nonlinear Schr\"odinger equation}, {\it Physica D}, {\bf 35} (1989), 167--185.

\bibitem{F1997}
A. S. Fokas, A unified transform method for solving linear and certain nonlinear PDEs, 
{\it Proc. Roy. Soc. Lond.} A {\bf 453} (1997), 1411--1443.

\bibitem{F2000}
A. S. Fokas, On the integrability of linear and nonlinear partial differential equations, {\it J. Math. Phys.} {\bf 41} (2000), 4188--4237.

\bibitem{F2002}
A. S. Fokas, Integrable nonlinear evolution equations on the half-line, 
{\it Comm. Math. Phys.} {\bf 230} (2002), 1--39.

\bibitem{F2002IMA}
A. S. Fokas, A new transform method for evolution PDEs, {\it IMA J. Appl. Math.} {\bf 67} (2002), 1--32.

\bibitem{F2003}
A. S. Fokas, The Davey-Stewartson I equation on the quarter plane with homogeneous Dirichlet boundary conditions, Integrability, topological solitons and beyond, {\it J. Math. Phys.} {\bf 44} (2003), 3226--3244.

\bibitem{F2004}
A. S. Fokas, Linearizable initial boundary value problems for the sine-Gordon equation on the half-line, {\it Nonlinearity} {\bf 17} (2004), 1521--1534.

\bibitem{F2005}
A. S. Fokas, A generalised Dirichlet to Neumann map for certain nonlinear evolution PDEs, 
{\it Comm. Pure Appl. Math.} {\bf LVIII} (2005), 639--670.

\bibitem{F2007}
A. S. Fokas, From Green to Lax via Fourier, {\it Proc. Sympos. Appl. Math.} {\bf 65}, Amer. Math. Soc., Providence, RI, 2007. 

\bibitem{F2009}
A. S. Fokas, The Davey-Stewartson on the half-plane, 
{\it Comm. Math. Phys.} {\bf 289} (2009), 957--993. 

\bibitem{Fbook}
A. S. Fokas, A unified approach to boundary value problems, CBMS-NSF regional conference series in applied mathematics, SIAM (2008).

\bibitem{FG1994}
A. S. Fokas and I. M. Gelfand, Integrability of linear and nonlinear evolution equations and the associated nonlinear Fourier transforms, {\it Lett. Math. Phys.} {\bf 32} (1994), 189--210.
 
\bibitem{FI1992} 
A. S. Fokas and A. R. Its, An Initial-Boundary Value Problem for the sine-Gordon Equation in Laboratory Coordinates, {\it Theor. Math. Phys.} \textbf{92} (1992), 387--403.
 
\bibitem{FI1992b} 
A. S. Fokas and A. R. Its,  Soliton generation for initial-boundary value problems, {\it Phys. Rev. Lett.} {\bf 68} (1992), 3117--3120. 

 \bibitem{FI1994}
A. S. Fokas and A. R. Its, An initial-boundary value problem for the Korteweg-de Vries equation, 
{\it Math. Comput. Simulation} {\bf 37} (1994), 293--321.

\bibitem{FI1996}
A. S. Fokas and A. R. Its, The linearization of the initial-boundary value problem of the nonlinear Schr\"odinger equation,
{\it SIAM J. Math. Anal.} {\bf 27} (1996), 738--764. 

\bibitem{FI2004}
A. S. Fokas and  A. R. Its, The nonlinear Schr\"odinger equation on the interval, {\it J. Phys. A.} {\bf 37} (2004), 6091--6114.

\bibitem{FIS2005}
A. S. Fokas, A. R. Its, and L.-Y. Sung, The nonlinear Schr\"odinger equation on the half-line, 
{\it Nonlinearity} {\bf 18} (2005), 1771--1822.

\bibitem{FK2004}
A. S. Fokas and S. Kamvissis, Zero-dispersion limit for integrable equations on the half-line with linearisable data, {\it Abstract Appl. Anal.} {\bf 5} (2004), 361--370.

\bibitem{FK2000}
A. S. Fokas and A. A. Kapaev, A Riemann-Hilbert approach to the Laplace equation, {\it J. Math. Anal. Appl.} {\bf 251} (2000), 770--804.
 
\bibitem{FK2003}
A. S. Fokas and A. A. Kapaev, On a transform method for the Laplace equation in a polygon, {\it IMA J. Appl. Math.} {\bf 68} (2003), 355--408.
 
\bibitem{FL2010}
A. S. Fokas and J. Lenells, Explicit soliton asymptotics for the Korteweg-de Vries equation on the half-line, {\it Nonlinearity} {\bf 23} (2010), 937--976.

\bibitem{FP2005}
A. S. Fokas and B. Pelloni, A transform method for linear evolution PDEs on a finite interval, {\it IMA J. Appl. Math.} {\bf 70} (2005), 564--587.

\bibitem{FPpreprint} 
A. S. Fokas and B. Pelloni, The Dirichlet-to-Neumann map for the elliptic sine-Gordon equation, submitted.

\bibitem{FPLpreprint} 
A. S. Fokas, B. Pelloni, and J. Lenells, Boundary value problems for the elliptic sine-Gordon equation in a semi-strip, submitted.

\bibitem{FP2006}
A. S. Fokas and D. A. Pinotsis, The Dbar formalism for certain linear non-homogeneous elliptic PDEs in two dimensions, {\it European J. Appl. Math.} {\bf 17} (2006), 323--346.

\bibitem{H1990}
I. T. Habibullin, B\"acklund transformations and integrable initial-boundary value problems, Nonlinear and Turbulent Processes vol 1 (Singapore: World Scientific) 1990, pp 130--138.

\bibitem{HZ2011}
M. H. Huang and Y. Q. Zhao, High-frequency asymptotics for the modified Helmholtz equation in a quarter-plane, {\it Appl. Anal.}, doi:10.1080/00036811.2010.534858.

\bibitem{I2003}
A. R. Its, The Riemann-Hilbert problem and integrable systems, {\it Notices Amer. Math. Soc.} {\bf 50} (2003), 1389--1400. 

\bibitem{K2010}
K. Kalimeris, Explicit soliton asymptotics for the nonlinear Schr\"odinger equation on the half-line, {\it J. Nonlinear Math. Phys.} {\bf 17} (2010), 445--452.

\bibitem{KF2010}
K. Kalimeris and A. S. Fokas, The heat equation in the interior of an equilateral triangle, {\it Stud. Appl. Math.} {\bf 124} (2010), 283--305.

\bibitem{K2003}
S. Kamvissis, Semiclassical nonlinear Schr\"odinger on the half line,
{\it J. Math. Phys.} {\bf 44} (2003), 5849--5868. 

\bibitem{Lax1968}
P. Lax, Integrals of nonlinear equations of evolution and solitary waves, {\it Comm. Pure Applied Math.} {\bf 21} (1968), 467--490.

\bibitem{Ldnls}
J. Lenells, The derivative nonlinear Schr\"odinger equation on the half-line, {\it Phys. D} {\bf 237} (2008), 3008--3019. 

\bibitem{LdnlsD2N}
J. Lenells, The solution of the global relation for the derivative nonlinear Schr\"odinger equation on the half-line, {\it Phys. D} {\bf 240} (2011), 512--525.

\bibitem{Lholedisk}
J. Lenells, Boundary value problems for the stationary axisymmetric Einstein equations: a disk rotating around a black hole, {\it Comm. Math. Phys.} {\bf 304} (2011), 585--635.

\bibitem{LFgnls}
J. Lenells and A. S. Fokas, An integrable generalization of the nonlinear Schr\"odinger equation on the half-line and solitons, {\it Inverse Problems} {\bf 25} 115006, 32pp.

\bibitem{LFernst}
J. Lenells and A. S. Fokas, Boundary-value problems for the stationary axisymmetric Einstein equations: a rotating disc, {\it Nonlinearity} {\bf 24} (2011), 177--206.

\bibitem{LFtperiodic}
J. Lenells and A. S. Fokas, The unified method: II NLS on the half-line with $t$-periodic boundary conditions, preprint.

\bibitem{LFinterval}
J. Lenells and A. S. Fokas, The unified method: III Non-linearizable problems on the interval, preprint.

\bibitem{MF2011}
D. Mantzavinos and A. S. Fokas, The Kadomtsev-Petviashvili II equation on the half-plane, {\it Phys. D} {\bf 240} (2011), 477--511. 

\bibitem{P2011}
V. G. Papanicolaou, An example where separation of variables fails, {\it J. Math. Anal. Appl.} {\bf 373} (2011), 739--744. 

\bibitem{PK2009}
T. S. Papatheodorou and A. N. Kandili,
Novel numerical techniques based on Fokas transforms, for the solution of initial boundary value problems, {\it J. Comput. Appl. Math.} {\bf 227} (2009), 75--82.

\bibitem{P2004}   
B. Pelloni, Well posed boundary value problems for linear evolution equations in finite intervals, {\it Math. Proc. Camb. Phil. Soc.} {\bf 136} (2004), 361--382.
    
\bibitem{P2005}
B. Pelloni, The spectral representation of two-point boundary-value problems for third-order linear evolution partial differential equations, {\it Proc. R. Soc. Lond. Ser. A} {\bf 461} (2005), 2965--2984.
   
\bibitem{P2006} 
B. Pelloni, Linear and nonlinear generalized Fourier transforms, {\it Philos. Trans. R. Soc. Lond. Ser. A} {\bf 364} (2006), 3231--3249.
     
\bibitem{PP2010}
B. Pelloni and D. A. Pinotsis, The elliptic sine-Gordon equation in a half plane, 
{\it Nonlinearity} {\bf 23} (2010), 77--88.
    
\bibitem{P2007}
D. A. Pinotsis, The Riemann-Hilbert formalism for certain linear and nonlinear integrable PDEs, {\it J. Nonlinear Math. Phys.} {\bf 14} (2007), 466--485.

\bibitem{S2005}
P. C. Sabatier, Generalised inverse scattering for a linear PDE associate to KdV, {\it J. Nonlinear Math. Phys.} {\bf 12} (2005), 599--613.

\bibitem{S2006}
P. C. Sabatier, Generalized inverse scattering transform applied to linear partial differential equations, {\it Inverse Problems} {\bf 22} (2006), 209--228.
   
\bibitem{S2011}
D. A. Smith, Spectral theory of ordinary and partial linear differential operators on finite intervals, PhD thesis, University of Reading, 2011.

\bibitem{Spreprint}
D. A. Smith, Well-posed two-point initial-boundary value problems with arbitrary boundary conditions, preprint.
   
\bibitem{SF2010a} 
E. A. Spence and A. S. Fokas, A new transform method I: Domain dependent
fundamental solutions and integral representations, 
{\it Proc. R. Soc. Lond. Ser. A} {\bf 466} (2010), 2259--2281.

\bibitem{SF2010b}
E. A. Spence, A. S. Fokas, A new transform method II: The global
relation and boundary value problems in polar co-ordinates, 
{\it Proc. R. Soc. Lond. Ser. A} {\bf 466} (2010), 2283--2307.


\bibitem{T1988}
V. O. Tarasov, An initial-boundary value problem for the nonlinear Schr\"odinger equation, {\it Zap. Nauchn. Sem. LOMI} {\bf 169} (1988), 151--165.

\bibitem{TF2008}
P. A. Treharne and A. S. Fokas, The generalized Dirichlet to Neumann map for the KdV equation on the half-line,
{\it J. Nonlinear Sci.} {\bf 18} (2008), 191--217. 

\bibitem{V2007}
S. Vetra, The computation of spectral representations for evolution PDE, PhD Thesis, University of Reading, 2007.

\end{thebibliography}

\end{document}